# Cumulative Flow Diagram-Based Fixed-Time Signal Timing Optimization at Isolated Intersections Using Connected Vehicle Trajectory Data

Chaopeng Tan, Yumin Cao, Xuegang (Jeff) Ban, Keshuang Tang*

***Abstract*—Time-dependent fixed-time control is a cost-effective control method that is widely employed at signalized intersections in numerous countries. Existing optimization models rely on traditional delay models with specific assumptions regarding vehicle arrivals. Recent advancements in intelligent mobility have led to development of high-resolution trajectory data of connected vehicles (CVs), thereby providing opportunities for improving fixed-time signal control. Taking advantage of CV trajectories, this study proposes a cumulative flow diagram (CFD)-based signal timing optimization method for fixed-time signal control at isolated intersections, which includes a CFD model and a multi-objective optimization model. The CFD model is formulated to profile the time-dependent vehicle arrival and departure processes under varying signal timing plans, where the intersection demand is estimated based on a weighted maximum likelihood estimation method. Then, a CFD-based multi-objective optimization model is proposed for both undersaturated and oversaturated traffic conditions. The primary objective is to minimize the exceeded queue dissipation time, whereas the secondary objective is to minimize the average delay at the intersection. Considering the data-driven property of the CFD model, a bi-level particle swarm optimization-based algorithm is then specially designed to solve the optimal cycle length (and the reference point if it is considered) and green ratios separately. The proposed method is evaluated based on simulation data and compared with Synchro. The results indicate that the proposed method outperforms Synchro under various traffic conditions in terms of average delay and queue since the time-dependent vehicle arrivals during the cycle are considered in a CV trajectory data-driven way.**

## I. Introduction

SIGNALIZED intersections are one of the most common types of bottlenecks in urban road networks; thus, traffic signal control plays an important role in improving urban traffic flow efficiency. Traffic signal control is commonly classified into three categories: fixed-time control, actuated control, and adaptive control. The latter two can adjust traffic signal timings based on real-time traffic flow information collected via infrastructure-based sensors [1]. However, the real-world applications of the two control methods are often constrained by the limited coverage, malfunction, and high maintenance cost of the sensors [2, 3]. In contrast, fixed-time signal control has proven to be a cost-effective alternative and is widely implemented at signalized intersections in many Asian countries, such as Japan, Korea, and China [4]. Therefore, it is important to improve the performance of fixed-time signal control using the available technologies [5, 6].

Existing fixed-time signal control optimization methods commonly consider traffic volume as basic input and rely on traditional delay models, such as the Webster delay model and the Highway Capacity Manual delay model, as well as specific assumptions on vehicle arrivals, such as uniform and random arrivals [7]. In real-world applications, the collection of traffic volumes incurs a significant cost, traffic volume data are not always available [8], vehicle arrival assumptions often deviate from reality, and the traditional delay models require careful validation and generate unreasonable results, particularly for near-saturated and oversaturated conditions [9]. The above limitations significantly hinder the practical application of existing fixed-time signal control optimization methods.

In the past decade, there has been a significant advancement in intelligent mobility and connected vehicle technologies, enabling access to massive amounts of high-resolution vehicle trajectory data from mobile sensors, such as smartphones, and other global positioning system-equipped devices. In general, vehicles equipped with mobile sensors can transmit their instantaneous positions and speeds with a high frequency (i.e., every 1–5 s), referred to as connected vehicles (CVs). In comparison with traditional fixed-location sensors, CV trajectory data can serve as a cost-effective data source for monitoring network-wide traffic flows, thus demonstrating the significant potential for facilitating the dynamic evaluation and optimization of traffic signals [10, 11].

Numerous studies have been devoted to applying CV trajectory data for the estimation of traffic flow parameters and performance measures at signalized intersections, including queue length [12-15], traffic volume [16-19], and delay [20, 21]. These studies have reported promising results, demonstrating the considerable potential of CV trajectory data. Nevertheless,

Manuscript received December 11, 2020. This work was supported by the National Natural Science Foundation of China under Grant 61673302 and Grant U1764261. The work of Chaopeng Tan was supported by China Scholarship Council, during the visit of the University of Washington from December 2019 to December 2021, under Grant 201906260111. (*Corresponding author: Keshuang Tang*)

Chaopeng Tan, Yumin Cao, and Keshuang Tang are with the College of Transportation Engineering, Tongji University, Shanghai 201804, China (e-mail: tanchaopeng@tongji.edu.cn; cao97@tongji.edu.cn; tang@tongji.edu.cn).

Xuegang Ban is with the Department of Civil and Environmental Engineering, University of Washington, Seattle, WA 98195 USA (e-mail: banx@uw.edu).

the methods proposed in these studies cannot be applied directly to signal control optimization because they do not provide analytical relationships between performance measures and varying signal timing parameters.

Some studies have focused on the optimization of fixed-time signal control based on CV trajectory data. Early studies mainly utilized the travel information provided by the CV trajectory data to optimize the offset for arterial coordination signal control [22, 23]. Yao *et al.* [24] proposed a STREAM model to optimize all signal timing parameters of arterial control, namely cycle length, green splits, and offset. This method first estimated the arrival rate between CVs, then achieved the evolution of the delay of CVs with signal timing parameters, and finally achieved arterial coordination control by minimizing the average delay of CVs. Aiming at network signal timing coordination, Yan *et al.* [25] proposed a network-level multiband signal coordination model. This study first used historical CV trajectories to identify major traffic streams in urban road networks and then optimized the progression bands for the identified major traffic streams.

Ma *et al.* [5] proposed a signal timing optimization model for isolated intersections. This study assumed that the vehicle arrival distribution is uniform, and CV trajectories during the cycle follow the same-ratio principles, which is the basis for the evolution of CV trajectories under varying signal timing parameters. However, vehicle arrivals are largely affected by signal control at the upstream intersections; thus, the assumption of uniform vehicle arrivals is ideal for practical applications. In addition, when varying the cycle length and green splits, the same-ratio principles of CVs can barely be satisfied because of varying vehicle arrivals during the cycle. In summary, the existing studies on fixed-time signal control at isolated intersections, either based on fixed-location detectors or CV trajectories, require specific assumptions on vehicle arrival or the penetration rate of CVs, which significantly constrains their real-world applications.

Meanwhile, a few studies have developed real-time signal control methods based on CV trajectory data [26-30], targeting future scenarios when CVs become popular. Nevertheless, because the information of regular vehicles needs to be estimated based on the information provided by CVs, these real-time control methods typically require a high penetration rate of CVs. Acceptable performance can only be achieved when the penetration rate of CVs exceeds 30% [5, 31]. Furthermore, the CV trajectory data are currently owned by ridesharing companies, whereas traffic signals are operated by traffic management sectors, and it is expected to take years to be able to apply the vehicle-to-infrastructure technology in a large-scale network [32]. Therefore, the implementation of CV-based real-time signal control is likely to be difficult even in the near future. In contrast, fixed-time signal control does not require real-time CV trajectory data and has great demand in the real world; therefore, it is worthwhile to develop fixed-time signal control optimization methods by utilizing offline CV trajectory data.

Therefore, this study proposes a CV trajectory-driven cumulative flow diagram (CFD) model for fixed-time signal timing optimization at isolated intersections. CFD represents a traffic flow profile that illustrates all vehicle accumulation and dissipation processes of a signalized intersection. In our previous research, we performed CFD estimation and prediction for a single phase based on CV trajectories for the initial and candidate signal timing plans, which are nonlinear and non-parametric [33]. In this study, a general linearized CFD model formulating CFD evolution under varying signal timing plans is first proposed, where the traffic demand of the intersection is estimated based on a novel statistical method, namely weighted maximum likelihood estimation (WMLE). This CFD model is CV trajectory data-driven and requires no prior assumptions regarding vehicle arrival and the penetration rate of CVs. Then, based on the proposed CFD model, a multi-objective optimization model is proposed to optimize the signal timing plans for both undersaturated and oversaturated conditions. In the multi-objective optimization model, the primary objective is to ensure that the queues are clear before the yellow start time, and the secondary objective is to minimize the average delay of the intersection. To solve this multi-objective optimization problem, a bi-level particle swarm optimization (PSO)-based algorithm is designed to determine the optimal cycle length and its corresponding green splits.

The major contributions of this study are summarized as follows:

1) We propose a novel statistical method for arrival rate estimation that relaxes the first-in-first-out assumption on vehicle arrivals, and we demonstrate that by considering the number of CVs, the estimation errors of phases can be neutralized, and the overall estimation accuracy of the intersection can be significantly improved.
2) Unlike existing delay models requiring pre-defined arrival patterns, we proposed a CV trajectory data-driven CFD model that can accurately profile vehicle arrival and departure processes under varying signal timing plans, which models time-dependent vehicle arrivals during the cycle without any prior assumptions on vehicle arrival patterns.
3) Based on the CFD model, we developed a hierarchical multi-objective signal timing optimization model for both undersaturated and oversaturated traffic conditions, in which both the primary and secondary objectives, i.e., exceeded queue dissipation time and average delay, are specially derived from the CFD model. Since the time-dependent vehicle arrivals during the cycle are considered, the proposed CFD-based optimization model can achieve extra signal control benefits by optimizing the signal timing parameters, particularly the reference point of the signal timing plan.
4) In comparison with Synchro, which uses ground truth volumes as input, the proposed CFD-based optimization model uses CV trajectories as input and can achieve better performances in terms of average delay and queue length under various traffic conditions, which provides a more cost-effective means for fixed-time signal timing optimization at isolated intersections.

## II. Cumulative Flow Diagram Model

In our previous study, we proposed a CV trajectory-driven CFD estimation and prediction framework for a single phase at signalized intersections [33]. A functional relationship between CFD and signal timing parameters was established; however, it is nonlinear and non-parametric. In addition, the previous framework requires the red start time of the phase to be set as the cycle start time, which implies that the cycle start time varies for each phase. These limitations make it challenging to apply the previous CFD framework to direct signal timing optimization. Thus, in this section, we improve the previous CFD framework and propose a general linearized and parameterized CFD model. It should be stated that, in this study, we will use a common CFD in the cycle to represent the traffic flow profile for the whole analysis period. That is to say, like other fixed-time signal timing optimization models, we assume that the vehicle arrivals for all cycles are the same. Nevertheless, the significance of this study is that we have considered the detailed time-dependent vehicle arrivals during the cycle.

The basic idea of the CFD model is as follows: First, by approximating the CV arrivals based on a piecewise constant function, the cumulative arrival curve $S_a$ is modeled as a piecewise linear function associated with the cycle length *C* and changes in the reference point $\Delta\varphi$ of the analysis period; then, the queue discharging wave speed is fitted by Huber regression, and the cumulative departure curve $S_d$ is modeled as a linear function of the green phase parameters, that is, effective green start time $t_{egs}$ and end time $t_{ege}$; eventually, the CFD model is obtained based on $S_a$ and $S_d$, and several performance measures, such as volume, queue length, total delay, and average delay, can be further derived. The flowchart of the proposed CFD model is shown in Fig. 1.

### *A. Cumulative arrival curve*

For urban signalized intersections, owing to the upstream signals, vehicle arrivals at the downstream intersection could be varied during different stages of the cycle. Hence, for a specific phase, a time-dependent arrival rate function $\lambda(t)$ of the cycle was adopted to represent the traffic arrival pattern during the analysis period,

$$\lambda(t) = \lambda_0 \cdot \alpha(t) \qquad 0 < t \le C \tag{1}$$

where $\lambda_0$ denotes the average arrival flow rate of the subject phase, which will be estimated later based on CV trajectories under the initial signal timing plan; $\alpha(t)$ denotes the time-dependent scaling parameter of the arrival rate during the cycle, which is generated by CV trajectories; and $C$ denotes the candidate cycle length.

Since the integral of $\lambda(t)$ over the time of the cycle should be equal to the traffic volume of the cycle, i.e., $\lambda_0 C$, we have $\int_0^C \lambda(t)dt = \int_0^C \lambda_0 \alpha(t)dt = \lambda_0 C$ , then we can obtain $\int_0^C \alpha(t)dt = C$.

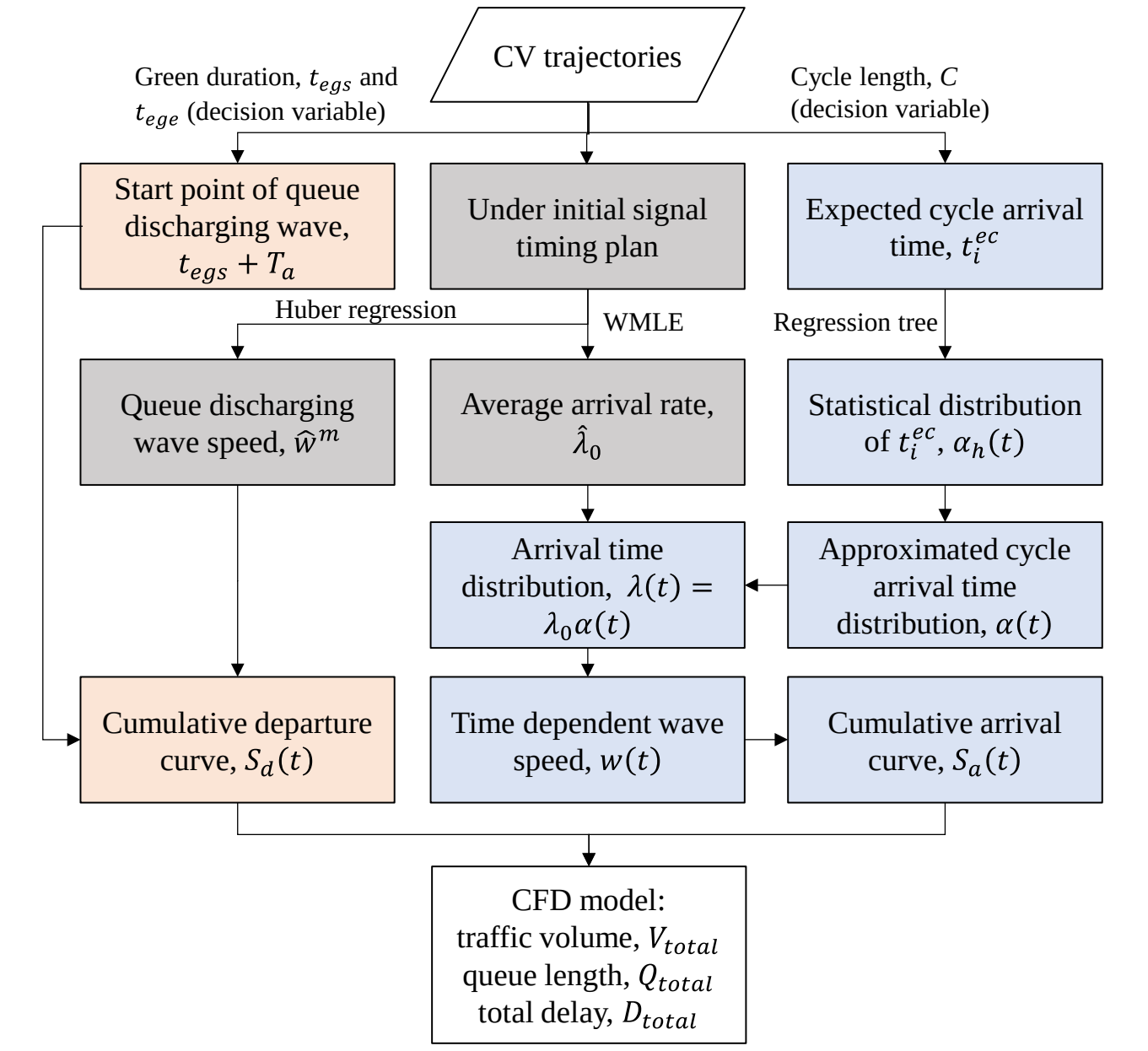

**Fig. 1.** Flowchart of the proposed CFD model.

Several studies have utilized sampled vehicle trajectories to obtain $\alpha(t)$ to estimate the average traffic volume [19] or cycle-based back of queues [12]. Typically, the arrival time distribution of CVs during the cycle is used to represent $\alpha(t)$ through non-parametric estimation methods such as histogram [19] or kernel density estimation [12, 33]. However, the non-parametric and nonlinear characteristics of $\alpha(t)$ make the cumulative arrival curve non-parametric and nonlinear as well. This causes difficulties in further signal timing optimization. Considering the simplicity of the solution and the applicability of the method, we used a piecewise constant function to fit $\alpha(t)$, transforming the cumulative arrival curve to be linear and parametric, thereby facilitating final signal timing optimization.

As shown in Fig. 2(a), for the trajectory of the *i*-th CV observed during the TOD period, given any cycle length *C*, its expected cycle arrival time $t_i^{ec}$ can be determined as follows,

$$t_i^{ec} = t_i^e - \Delta\varphi - \left\lfloor \frac{t_i^e - \Delta\varphi}{C} \right\rfloor \times C \tag{2}$$

where $t_i^e$ denotes the expected arrival time of the *i*-th CV observed during the TOD period; $\Delta\varphi$ denotes the change in the reference point, $0 \le \Delta\varphi < C$, and the reference point indicates the start time of the first cycle during the analysis period; $\lfloor\cdot\rfloor$ indicates the integer part.

Then, the statistical distribution of $t_i^{ec}$ of all CVs, that is, $\alpha_h(t)$, can be obtained, as shown by the histogram in Fig. 2(b). Nevertheless, the statistical distribution fluctuates significantly and is discontinuous and nonlinear. Therefore, we use a piecewise constant function to approximate the statistical distribution, as shown in Fig. 2(b). Note that when the penetration rates are low (e.g., no more than 10%), very few CVs would be observed during the analysis period for those phases with low volume. This makes $\alpha_h(t)$ sparser and more

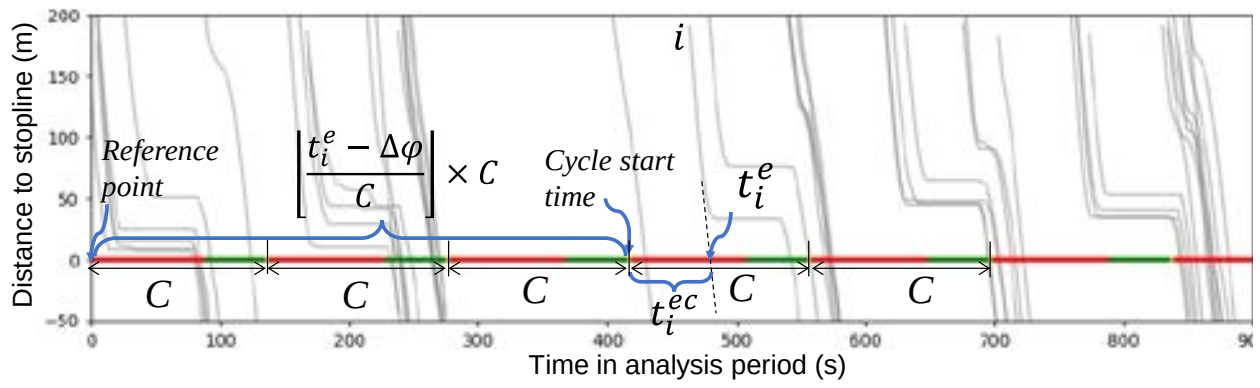

(a) Calculation of the expected cycle arrival time of CVs

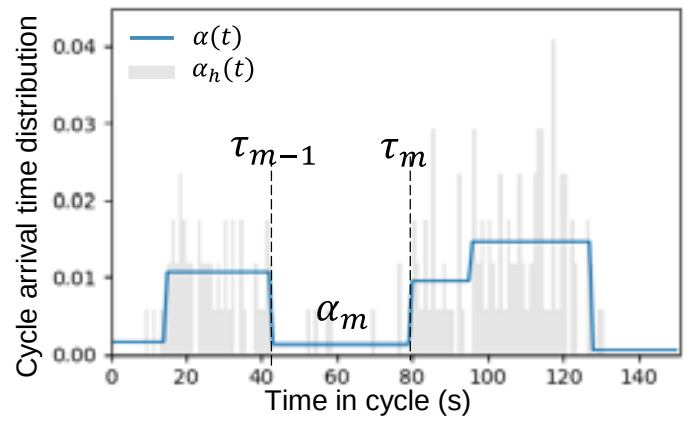

(b) Piecewise constant function of cycle arrival rates

**Fig. 2.** Approximation of the arrival time distribution.

challenging to fit directly. Thus, under such conditions, the application of a kernel density estimator is first recommended to smooth and fit the data [33, 34], and then the use of a piecewise constant function is recommended to approximate the KDE fitted distribution. Given any value of $M$, the piecewise constant function $\alpha(t)$ can be estimated by solving the following optimization problem with the objective of minimizing the difference between $\alpha(t)$ and $\alpha_h(t)$,

$$P1: \min_{\alpha_m, \tau_m} \sum_{m=1}^{M} \left[ \sum_{t=\tau_{m-1}}^{\tau_{m-1}} (\alpha_h(t) - \alpha_m)^2 \right] \quad (3)$$

where $\tau_0$ is equal to 0; $\tau_M$ is equal to the cycle length $C$; $M$ is the total number of segments of the piecewise constant function, whose optimal value is discussed in Section 4.2.1; and $\alpha_m$ is the value of $\alpha(t)$ during the period $[\tau_{m-1}, \tau_m)$.

The key to solving P1 is to determine the best split points, that is, $\tau_m$, $m = 1,2, \ldots, M-1$. Then, the corresponding value $\alpha_m$ for periods $[\tau_{m-1}, \tau_m), m = 1,2, \ldots, M$ can be easily determined by minimizing the sum of squared errors (SSE). However, it is computationally infeasible to consider every possible combination of the split points. Therefore, a tree-based method, the regression tree, is adopted to solve the optimal split points [35]. The splitting process is achieved by recursive binary splitting, which is a top-down greedy approach. In our case, this approach proceeds as follows:

1. First, given $\hat{\alpha}_h(t)$, we search for an optimal split time point $st_1$ that minimizes the overall SSE as

$$st_1 = \arg\min_{t_x} \left\{ \sum_{t=0}^{t_x} (\hat{\alpha}_h(t) - \bar{\alpha}_1)^2 + \sum_{t=t_x+1}^{C} (\hat{\alpha}_h(t) - \bar{\alpha}_2)^2 \right\} \quad (4)$$

where $\bar{\alpha}_1$ and $\bar{\alpha}_2$ are the average of two parts, respectively.

2. Then, we repeat this process on the part with the largest SSE until the stopping criterion associated with the overall SSE differences is reached or $M-1$ split points are obtained.
3. Eventually, after we partition $\hat{\alpha}_h(t)$ into $M$ parts in terms of time, the estimator of each part is determined by minimizing the SSE.

Therefore, a piecewise constant function of $\alpha(t)$ in the general form can be obtained as follows,

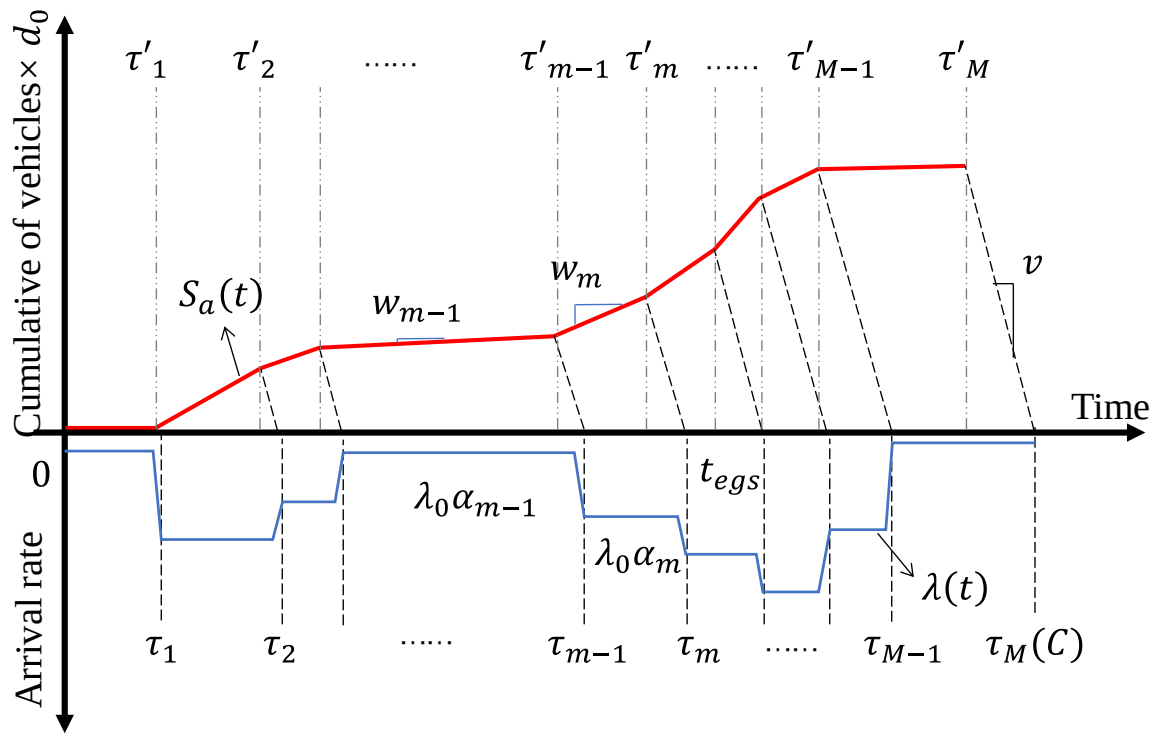

**Fig. 3.** Derivation of the cumulative arrival curve.

$$\alpha(t) = \begin{cases} \alpha_1, & \tau_0 \le t < \tau_1 \\ \alpha_2, & \tau_1 \le t < \tau_2 \\ & \ldots \\ \alpha_M, & \tau_{M-1} \le t < \tau_M. \end{cases} \quad (5)$$

As we derived that $\int_0^C \alpha(t)dt = C$, we have $\sum_{m=1}^{M} [\alpha_m(\tau_m - \tau_{m-1})] = C$.

In summary, once we estimate $\lambda_0$ based on CV trajectories under the initial signal timing plan and generate $\alpha(t)$ by approximating the arrival time distribution of CVs, $\lambda(t)$ for any candidate $C$ and $\Delta\varphi$ can be obtained.

In Tan *et al.* [33], the time mapping relation between the instantaneous wave speed and arrival rate was derived, as presented in Eq. (6). For a specific moment $t$ in the cycle, whose instantaneous arrival rate is $\lambda(t)$, its corresponding wave speed should be $w(t')$ at moment $t'$. The detailed derivation process can be found in Tan *et al.* [33].

$$w(t') = \frac{v^2}{v - d_0\lambda(t)} - v, \quad 0 \le t \le C \quad (6)$$

where $d_0$ denotes the vehicle jam spacing; $v$ denotes the free-flow speed; $t' = t(1 - d_0\bar{\lambda}_t/v)$ ; $t' \in [0, C(1 - d_0\lambda_0/v)]$; and $\bar{\lambda}_t = \lambda_0 \int_0^t \hat{\alpha}(s)ds/t$.

Thus, given the piecewise constant $\alpha(t)$, the wave speed is also a piecewise constant function,

$$w(t) = \begin{cases} \frac{vd_0\lambda_0\alpha_1}{v - d_0\lambda_0\alpha_1}, & \tau'_0 \le t < \tau'_1 \\ \frac{vd_0\lambda_0\alpha_1}{v - d_0\lambda_0\alpha_2}, & \tau'_1 \le t < \tau'_2 \\ \ldots \\ \frac{vd_0\lambda_0\alpha_M}{v - d_0\lambda_0\alpha_M}, & \tau'_{M-1} \le t < \tau'_M \end{cases} \quad (7)$$

where $\tau'_m = \tau_m - d_0\lambda_0 \sum_{\epsilon=1}^{m} [\alpha_\epsilon(\tau_\epsilon - \tau_{\epsilon-1})]/v, m = 1,2, \ldots, M$.

Finally, as shown in Fig. 3, by integrating the wave speed in terms of time, a piecewise linear cumulative arrival curve $S_a$ of vehicles can be acquired as

$$S_a(t) = w_m t - w_m \tau'_{m-1} + \sum_{\epsilon=1}^{m-1} w_\epsilon(\tau'_\epsilon - \tau'_{\epsilon-1}),$$
$$\tau'_{m-1} \le t < \tau'_m, m = 1,2, \ldots, M \quad (8)$$

In particular, $m-1$ should be no less than 1 for $\sum_{\epsilon=1}^{m-1} X$; however, in this study, we assume that when $m-1=0$, we have $\sum_{\epsilon=1}^{0} X = 0$ for any $X$. Thus, when $m=1$, Eq. (8) is true.

*B. Cumulative departure curve*

To model the cumulative departure curve under candidate signal plans, we must first estimate the queue discharging wave speed, that is, the slope of the cumulative departure curve $w^m$, under the initial signal timing plan. Tan *et al.* [33] proposed a robust regression method, that is, Huber regression [36], for queue discharging wave speed estimation. The evaluation results indicate that Huber regression is robust against outliers (such as abnormal discharging points deviating significantly from the majority and those caused by other traffic events near an intersection) and can significantly improve the accuracy of CFD estimation. Similar to Tan *et al.* [33], the Huber regression is also used to estimate $w^m$ under the initial signal timing plan, and we can obtain the following:

$$S_d^0(t) = \widehat{w}^m t + \hat{\xi}^0 \tag{9}$$

where $S_d^0$ represents the cumulative departure curve under the initial signal timing plan; $\widehat{w}^m$ denotes the queue discharging wave speed; and $\hat{\xi}^0$ is a constant value.

Notably, for a specific cycle, the green start time of the subsequent phase is determined by the green time of the previous phase. In other words, if we consider the red start time of each phase as the origin of the corresponding cumulative arrival curve, the cycle origin of each phase differs. Then, for each feasible solution of cycle length and green splits, $\alpha(t)$ must be regenerated for all phases, which significantly increases the computational burden and makes the signal timing optimization problem challenging to solve. Therefore, we set the green start time of the first phase in the cycle as the common origin for all phases. In this manner, given a feasible $C$, the $\alpha(t)$ for all phases only needs to be generated once. Thus, the cumulative departure curve $S_d$ for any studied phase can be modeled as a function of the corresponding effective green start time $t_{egs}$ and end time $t_{ege}$.

As depicted in Fig. 4, considering the general case in which the green start time of the studied phase is not the origin of the cycle, the red time is divided into two parts by the green time, that is, $R_a$ (or $C - t_{ege}$) and $R_b$ (or $t_{egs}$). Notably, the cumulative departure curve is not expected to start from $t_{egs}$ because the queue $Q_a$ accumulated during $R_a$ dissipates early in the green time.

$Q_a$ is calculated based on the number of vehicles arriving during $R_a$:

$$Q_a = S_a(C) - S_a(t_{ege})$$
$$= d_0\lambda_0 C - d_0\left[\lambda_0 \sum_{\epsilon=1}^{s} \alpha_\epsilon(\tau_\epsilon - \tau_{\epsilon-1}) + \lambda_0\alpha_{s+1}(t_{ege} - \tau_s)\right],$$
$$\tau_s \le t_{ege} < \tau_{s+1}, s = 0,1,\ldots,M-1 \tag{10}$$

where $\lambda_0 C$ denotes the total volume of the cycle, and the terms within square brackets represent the number of vehicles arriving before $t_{ege}$.

The corresponding dissipation time of $Q_a$ can be calculated as follows:

$$T_a = Q_a\left(\frac{1}{\widehat{w}^m} + \frac{1}{v_l}\right) \tag{11}$$

where $v_l$ denotes the average queue departure speed of queued vehicles and $v_l \le v$.

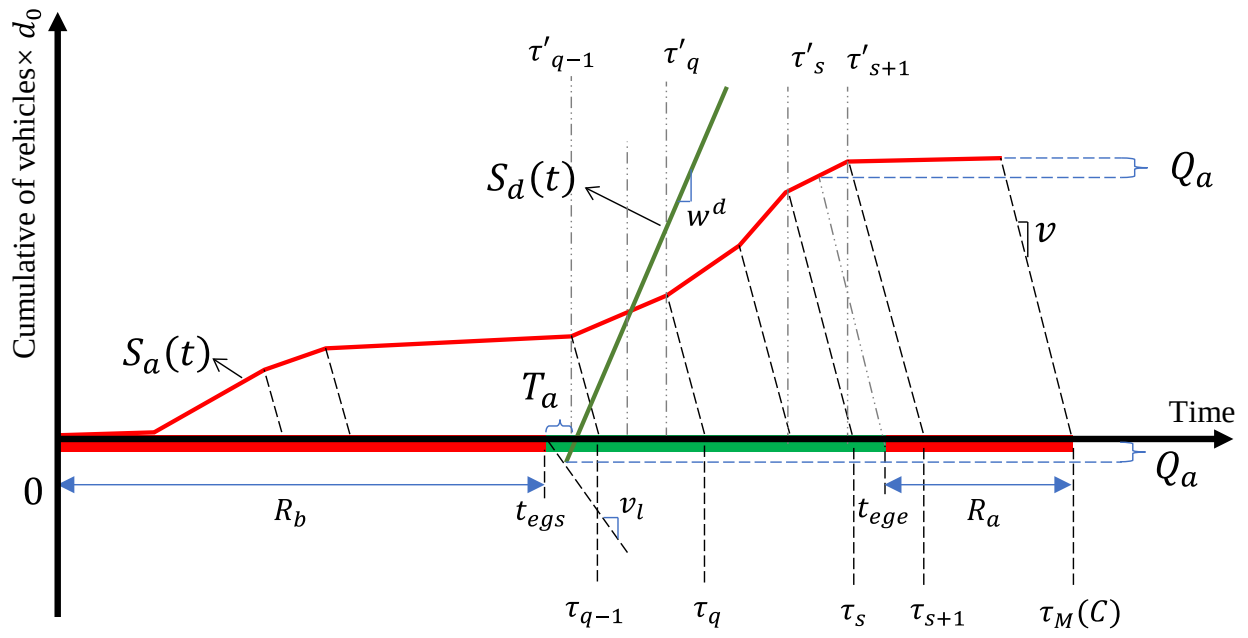

**Fig. 4.** Derivation of the cumulative departure curve.

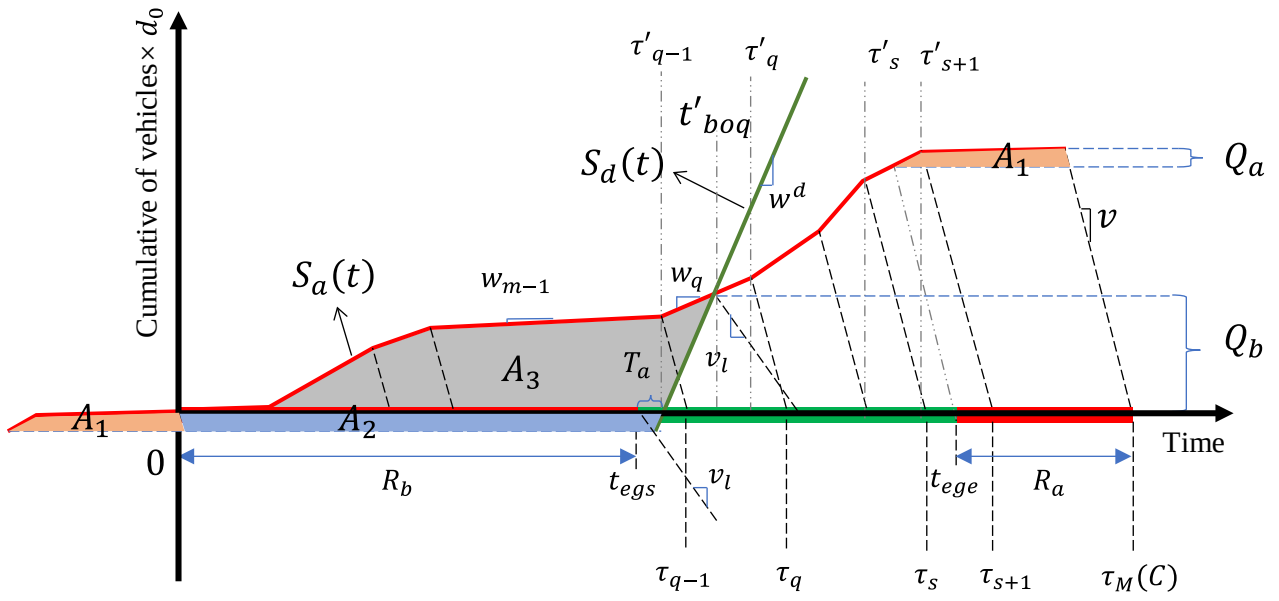

**Fig. 5.** Illustration of the total delay calculation with common cycle origin.

Thus, the cumulative departure curve can be determined as follows:

$$S_d(t) = \widehat{w}^m(t - t_{egs} - T_a) \tag{12}$$

*C. CFD model*

As both the cumulative arrival and departure curves were modeled as functions of the signal timing parameters, the CFD model was formulated. Based on this model, several critical performance measures for signal timing optimization can be derived as functions of the signal timing parameters.

As shown in Fig. 5, the queue length consists of two parts, that is, $Q_a$ and $Q_b$, where $Q_b$ is determined by the point at which $S_a(t)$ and $S_d(t)$ intersect, as follows:

$$t'_{boq} = \frac{\widehat{w}^m(t_{egs}+T_a) - w_q\tau'_{q-1} + S_a(\tau'_{q-1})}{\widehat{w}^m - w_q} \tag{13}$$

$$Q_b = \frac{w_q(t_{egs}+T_a) - w_q\tau'_{q-1} + S_a(\tau'_{q-1})}{1 - w_q/\widehat{w}^m} \tag{14}$$

where $t'_{boq}$ is assumed to be located in $[\tau'_{q-1}, \tau'_q), 1 \le q \le M$; $w_q$ is the wave speed between $[\tau'_{q-1}, \tau'_q)$; and $S_a(\tau'_{q-1}) = \sum_{\epsilon=1}^{q-1} w_\epsilon(\tau'_\epsilon - \tau'_{\epsilon-1}) = d_0\lambda_0 \sum_{\epsilon=1}^{q-1} \alpha_\epsilon(\tau_\epsilon - \tau_{\epsilon-1})$ as indicated by Eq. (8).

Note that the delay in this study refers to the stopped delay, which is the amount of time vehicles are actually stopped owing to the signal. The total delay consists of three parts, that is, $A_1, A_2$, and $A_3$, where $A_1 + A_2$ denotes the delay produced by vehicles in $Q_a$, and $A_3$ denotes the delay produced by vehicles in $Q_b$.

$$A_1 = 0.5 \sum_{\epsilon=s+1}^{M} \{(\tau_\epsilon - \tau_{\epsilon-1})[S_a(\tau'_\epsilon) + S_a(\tau'_{\epsilon-1}) - 2S_a(\tau'_s)]\}$$
$$-0.5(2C - \tau_s - t_{ege})[d_0\lambda_0\alpha_{s+1}(t_{ege} - \tau_s)] \tag{15}$$

$$A_2 = 0.5Q_a \times \left(2t_{egs} + T_a + \frac{Q_a}{v_l} - \frac{Q_a}{v}\right) \tag{16}$$

$$A_3 = 0.5\sum_{\epsilon=1}^{q-1}\{(\tau'_\epsilon - \tau'_{\epsilon-1})[S_a(\tau'_\epsilon) + S_a(\tau'_{\epsilon-1})]\} + 0.5(t'_{boq} - \tau'_{q-1})[S_a(\tau'_{q-1}) + Q_b] \quad (17)$$

In summary, we can obtain the following performance measures: traffic volume $V_{total}$, queue length $Q_{total}$, and total delay $D_{total}$:

$$V_{total}(C) = \lambda_0 C \quad (18)$$

$$Q_{total}(C, t_{egs}, t_{ege}) = Q_a + Q_b \quad (19)$$

$$D_{total}(C, t_{egs}, t_{ege}) = A_1 + A_2 + A_3 \quad (20)$$

*D. Parameter estimation for CFD*

In the CFD model, $\lambda_0$ is a critical parameter to be estimated because it reflects the traffic demands and, to some extent, determines the cycle length and green splits. Several studies have attempted to estimate the traffic volume based on CV trajectories [16-19, 33]. However, these studies still have limitations: 1) Existing studies can only estimate the traffic volume for undersaturated conditions and are not applicable for oversaturated conditions. 2) Existing statistical approaches, such as those proposed by Zheng and Liu [19] and Yao *et al.* [18], commonly adopt the relative arrival time and queuing position of consecutive CVs for model construction and assume that the queue process of vehicles conforms to the first-in-first-out (FIFO) principle, which is not always true in the case of multiple lanes. Data-driven methods, such as the one proposed by Tang *et al.* [16], are typically computationally burdensome. 3) Existing studies estimate the traffic volume for each phase separately. Some phases may be overestimated, whereas others may be underestimated. This may lead to a mismatch between the green ratios and the traffic demands of these phases during optimization, reducing the benefits of the optimized signal plan.

To overcome the aforementioned limitations, in this study, we introduced a novel statistical method, namely WMLE along with redistribution, to estimate the intersection arrival rates under both undersaturated and oversaturated conditions. First, the WMLE is adopted to initially estimate the arrival rates for each phase. When constructing the weighted maximum likelihood function, all CVs are treated as independent observations of the arrival rate based on their arrival time and queuing positions during the cycle, which relaxes the assumption of FIFO. Then, the number of CVs is adopted to redistribute the estimates of all phases, which can neutralize the estimation errors and improve the accuracy of the estimates of all phases.

As indicated by Eq. (1), the vehicle arrivals during the cycle are modeled using $\lambda_0\alpha(t)$. Given the initial signal timing plan, $\alpha(t)$ can be obtained using Eq. (5). Then, each queued CV can be treated as an independent observation of the arrival rate, as shown in Fig. 6. Thus, we have the probability density function of the *i-th* queued CV observed during the analysis period based on the time-dependent Poisson process [12, 15, 18, 19]:

$$f(n_i, t_i^{ec}|\lambda_0) = \frac{\left[\lambda_0\int_0^{t_i^{ec}}\alpha(t)dt\right]^{n_i}}{n_i!}e^{-\lambda_0\int_0^{t_i^{ec}}\alpha(t)dt} \quad (21)$$

where $n_i$ denotes the number of vehicles joining the queue before the *i-th* queued CV during the same cycle. For vehicles queued more than once, only the earliest queuing positions of these CVs were considered.

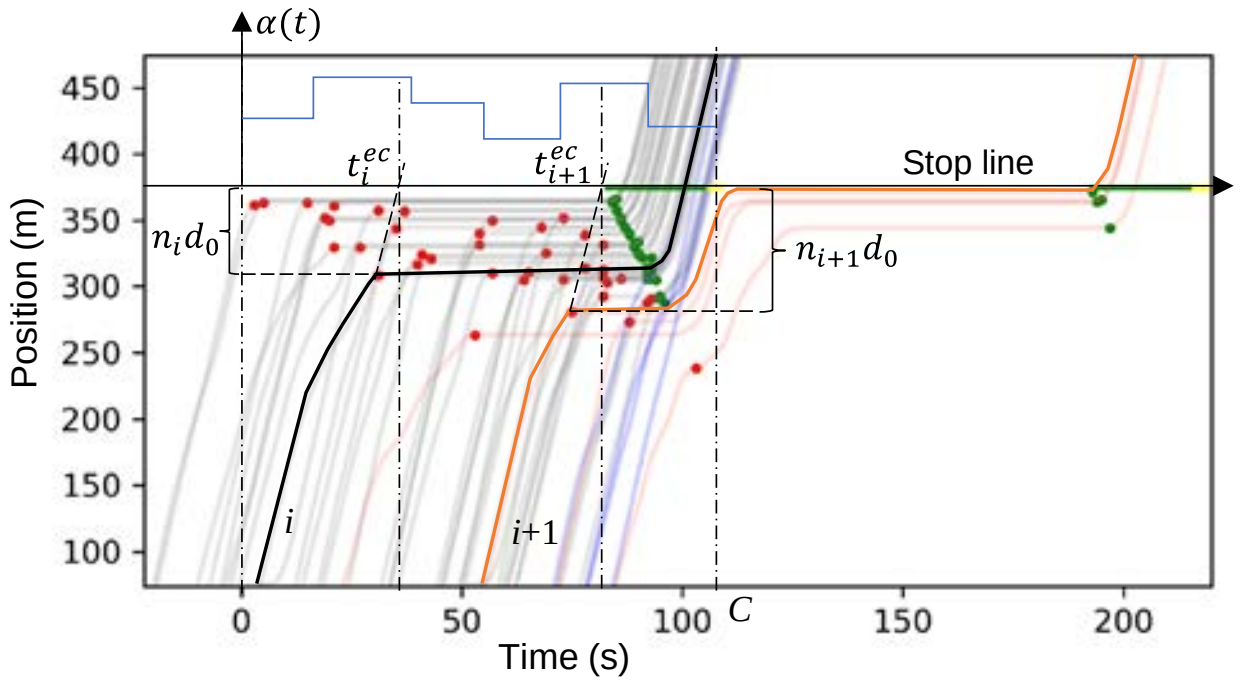

**Fig. 6.** Illustration of arrival rate estimation by WMLE.

Under constantly oversaturated conditions, the queue infinitely increases, and spillover may occur, making the proposed estimation method inapplicable. Under slightly or temporally oversaturated conditions, the overflow queue is relatively stable, although it still exists during the analysis period, and the proposed method is applicable. Under undersaturated conditions, $n_i$ denotes the exact queuing position of the CVs, while for temporally oversaturated conditions, the average overflow queue should be excluded:

$$n_i = \begin{cases} n_i^q, & undersaturated \\ n_i^q - n^{oq}, & oversaturated \end{cases} \quad (22)$$

where $n_i^q$ denotes the queuing position of the *i-th* queued CV and $n^{oq}$ denotes the average overflow queue of the analysis period. In this study, a queuing position is identified as an overflow queue position if more than 50% of the queued CVs observed at this position queue more than once, and $n^{oq}$ represents the maximum overflow queue position.

Then, the WMLE is adopted to formulate the weighted likelihood of $N_q$ queued CVs as [37]

$$L(\lambda_0) = \prod_{i=1}^{N_q} f(n_i, t_i^{ec}|\lambda_0)^{w_i}, \quad (23)$$

where $w_i$ denotes the weight of the *i-th* queued CV.

Note that a greater arrival time $t_i^{ec}$ covers a longer cycle length, indicating a more reliable observation of the arrival rate. In addition, the measurement error of $n_i$ has a low impact on the result with a greater $t_i^{ec}$. Therefore, in our case, considering the time-dependent characteristics of vehicle arrivals during the cycle, we use $\int_0^{t_i^{ec}}\alpha(t)dt$ as the weight of observations, and the log-likelihood can be written as

$$\log L(\lambda_0) = \sum_{i=1}^{N_q}\int_0^{t_i^{ec}}\alpha(t)dt \log f(n_i, t_i^{ec}|\lambda_0). \quad (24)$$

By computing the derivative of this log-likelihood, $\lambda_0$ was initially estimated as

$$\hat{\lambda}_0 = \frac{\sum_{i=1}^{N_q} n_i\int_0^{t_i^{ec}}\alpha(t)dt}{\sum_{i=1}^{N_q}\left(\int_0^{t_i^{ec}}\alpha(t)dt\right)^2} \quad (25)$$

Existing studies have shown that the lower the penetration rate, the more random the distribution of CVs in the population, which may result in more errors and variances in the estimation of $\lambda_0$. Similarly, the proposed $\lambda_0$ estimation method also faces the same problem. To optimize the signal timing plan, we must estimate $\lambda_0$ for each flow. Among the estimated arrival rates for different flows, some may be overestimated, while others may be underestimated. This may lead to a mismatch between the optimized green ratios and true traffic demands, reducing the benefits of signal control. To reduce the impact of this mismatch problem, we use the number of CVs to redistribute the estimated $\lambda_0$ for each phase as follows:

$$\tilde{\lambda}_0^k = \frac{1}{z^k} \times \left( \frac{0.5 z^k \hat{\lambda}_0^k}{\sum_{j=1}^{K} z^j \hat{\lambda}_0^j} + \frac{0.5 N^k}{\sum_{j=1}^{K} N^j} \right) \times \sum_{j=1}^{K} z^j \hat{\lambda}_0^j \tag{26}$$

where $\tilde{\lambda}_0^k$ denotes the redistributed arrival rate used for later signal timing optimization; $\hat{\lambda}_0^k$ is the arrival rate initially estimated by WML; $z^k$ is the lane number of phase $k$; and $N^k$ denotes the observed number of CVs of phase $k$; $K$ is the total number of phases in which the queued CVs were observed during the analysis period; $\sum_{j=1}^{K} z^j \hat{\lambda}_0^j$ indicates the total arrivals of $K$ phases; $\frac{z^k \hat{\lambda}_0^k}{\sum_{j=1}^{K} z^j \hat{\lambda}_0^j}$ indicates the proportion of the estimates of phase $k$ to the total; $\frac{N^k}{\sum_{j=1}^{K} N^j}$ indicates the proportion of the number of observed CVs of phase $k$ to the total; items in brackets indicate that we give the estimates the same weight as the number of trajectories to modify the initial estimates.

In corner cases, some phases may only have non-queued CVs due to well coordination between adjacent intersections, then the proposed WMLE method is not applicable. Their average arrival rates can be estimated by the number of CVs as follows:

$$\tilde{\lambda}_0^{k'} = \frac{1}{z^{k'}} \times \frac{N^{k'}}{\sum_{j=1}^{K} N^j} \times \sum_{j=1}^{K} z^j \hat{\lambda}_0^j \tag{27}$$

where $\tilde{\lambda}_0^{k'}$ is the average arrival rate of phase $k'$, which did not observe any queued CVs during the analysis period; $N^{k'}$ is the number of (non-queued) CVs of phase $k'$.

To sum up, as long as at least one phase has queue CVs, whose average arrival rate can be estimated by the WMLE method, we can still use Eq. (27) to obtain the estimates of the other phases without queued CVs. That is to say, our method fails only if all phases of the intersection have no queued CVs. Note that, in future real-world applications, the arrival rates estimated by other state-of-the-art methods or detected by other sensors can also be used as alternatives.

In addition to the average arrival rate, other basic parameters must be calibrated based on the observed CV trajectories. For each phase, the free-flow speed $v$ for each phase is adopted as the 95-th percentage speed value of all trajectory points of the corresponding phases; the average queue departure speed $v_l$ is adopted as the average speed value of the queued CVs when they travel through the stop line; and the saturated headway $h_s$ is calculated as follows:

$$h_s = d_0 \left( \frac{1}{\hat{w}^d} + \frac{1}{v_l} \right) \tag{28}$$

where vehicle jam spacing $d_0$ is empirically determined.

## III. CFD-Based Optimization Model

In this study, we applied the proposed method to a NEMA ring-barrier structure-controlled intersection, as illustrated in Fig. 7. The proposed optimization method can be extended to more general cases, but only the dual ring is discussed in this study as an example.

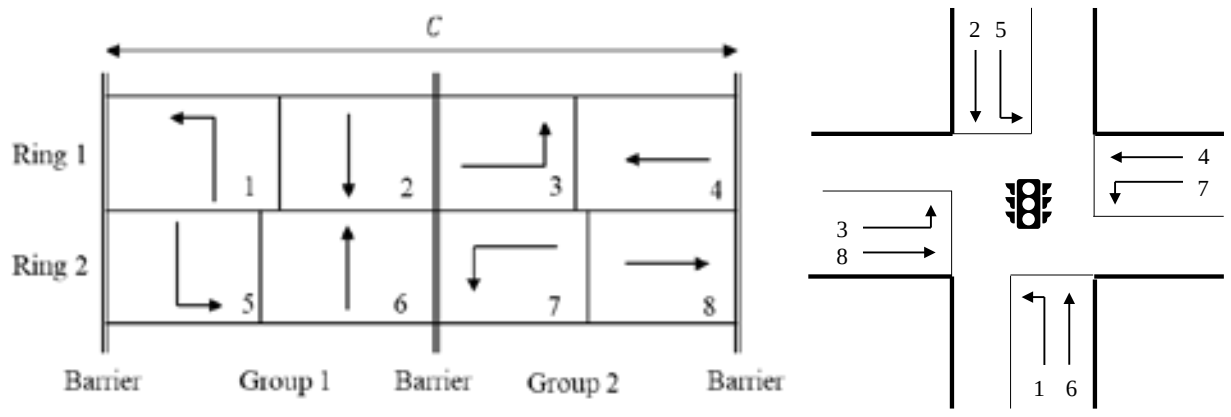

**Fig. 7.** Example phase structure.

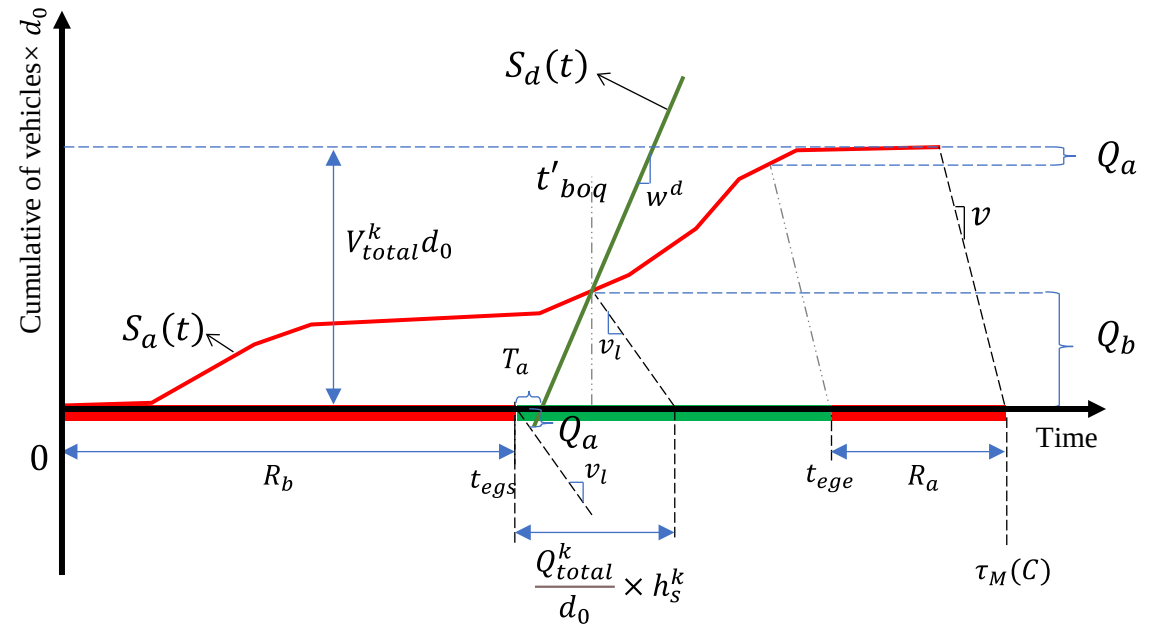

(a) Undersaturated

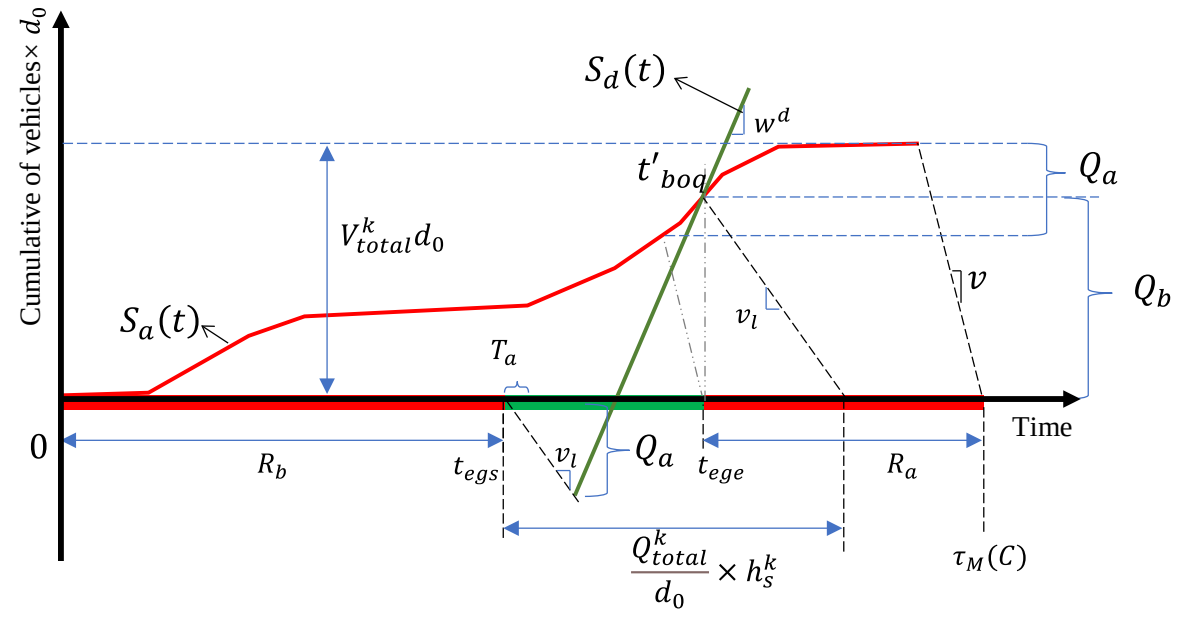

(b) Oversaturated

**Fig. 8.** Different traffic conditions for the CFD model.

### A. Objective

Given a feasible combination of cycle length and green ratios, each phase may have two states, that is, undersaturated or oversaturated, as shown in Fig. 8.

Typically, the existence of an overflow queue is used to identify whether the phase is oversaturated, as follows:

$$Q_{over}^k = \begin{cases} 0, & if\ \frac{Q_{total}^k h_s^k}{d_0} \le t_{ege}^k - t_{egs}^k \\ Q_{total}^k - \frac{(t_{ege}^k - t_{egs}^k) d_0}{h_s^k}, & others \end{cases} \tag{29}$$

If $Q_{over}^k > 0$, the phase $k$ is oversaturated; otherwise, it is undersaturated.

The proposed CFD model indicates the average traffic flow profile during the analysis period. The case of $Q_{over}^k > 0$ implies that the average demand exceeds the capacity, and theoretically, the queue increases infinitely. Under oversaturated conditions, the present CFD model cannot adequately represent the periodic traffic flow profile. However, although $Q_{over}^k$ does not represent a real overflow queue, it can still, to some extent, reflect the severity of the oversaturated condition of the phase. In contrast, when $Q_{over}^k = 0$, for the cases wherein the average demand is smaller than but close to the capacity, owing to the fluctuation of traffic demand, overflow queues may also exist for some cycles. Therefore, in this study, a more conservative criterion similar to $Q_{over}^{k}$ was adopted to determine the saturation condition of the studied phase, as follows:

$$G_{exceed}^{k} = \begin{cases} 0, \ if \ \frac{Q_{total}^{k} h_s^k}{d_0} + t_{egs}^{k} - g_e^k \le 0 \\ \frac{Q_{total}^{k} h_s^k}{d_0} + t_{egs}^{k} - g_e^k, \ others \end{cases} \tag{30}$$

where $G_{exceed}^{k}$ is defined as the exceeded queue dissipation time representing the difference between the queue dissipation time $Q_{total}^{k} h_s^k / d_0 + t_{egs}^{k}$ and the green end time $g_e^k$. If $G_{exceed}^{k} = 0$, the queue can dissipate before the green end time, and the phase is identified as an undersaturated condition; otherwise, it is oversaturated.

Thus, a hierarchical multi-objective structure is proposed for the signal timing optimization of both undersaturated and oversaturated traffic conditions:

$$\min_{C, \boldsymbol{g_s}, \boldsymbol{g_e}} f_{green} = \sum_{k=1}^{8} z^k G_{exceed}^{k} \tag{31a}$$

$$\min_{C, \boldsymbol{g_s}, \boldsymbol{g_e}} f_{delay} = \frac{\sum_{k=1}^{8} z^k D_{total}^{k}}{\sum_{k=1}^{8} z^k V_{total}^{k}} \tag{31b}$$

The primary objective is to ensure that the queues are clear before the yellow start time, which, to some extent, buffers the traffic fluctuation; this also ensures that the intersection is unsaturated by satisfying Eq. (28). The secondary objective is to minimize the theoretical average delay derived by the CFD model, which is commonly used for undersaturated conditions. The total delay and total volume derived by the CFD model are both values per lane per cycle; thus, they are weighted by the number of lanes.

Furthermore, we can scalarize the hierarchical multi-objective optimization problem into a single-objective optimization problem by weighting $f_{green}$ with a sufficiently large constant $\delta$, as follows:

$$P2: \min_{C, \boldsymbol{g_s}, \boldsymbol{g_e}} f_{hybrid} = f_{delay} + \delta f_{green} \tag{32}$$

## *B. Constraints*

Based on the example phase structure in Fig. 7, the following structural constraints can be formulated as follows:

$$\sum_{k=1}^{4} (g_e^k - g_s^k) = \sum_{k=5}^{8} (g_e^k - g_s^k) = C - 4 \times (y^k + r_a^k) \tag{33a}$$

$$g_e^1 - g_s^1 + g_e^2 - g_s^2 = g_e^5 - g_s^5 + g_e^6 - g_s^6 \tag{33b}$$

$$g_e^k + y^k + r_a^k = g_s^{k+1} \quad k = 1,2,3,5,6,7 \tag{33c}$$

$$g_e^k - g_s^k \ge g_{min}^k \quad k = 1,2,\ldots,8 \tag{33d}$$

$$C_{max} \ge C \ge C_{min} \tag{33e}$$

Eq. (33a), (33b), and (33c) confirm that the green phases should satisfy the structure of the dual ring, as shown in Fig. 7. In addition, the green time should be limited by a lower bound, which is determined by the time at which the pedestrian passes through the intersection when the pedestrian is considered, as shown in Eq. (33d). In general, the cycle length should be bounded within an appropriate range, as shown in Eq. (33e).

## *C. Solution*

As mentioned earlier, for each feasible cycle length, $\alpha(t)$ needs to be regenerated based on CV trajectories by solving P1; thus, the cumulative arrival curve varies for different cycle lengths owing to its data-driven property. This implies that it is impossible to calculate the gradient of the objective function, and it is quite difficult to solve P2 using exact algorithms. Therefore, an advanced meta-heuristic, PSO algorithm, is tailored into a bi-level form to solve P2, which solves the optimal cycle length and green splits separately, as shown in Fig. 9. The PSO algorithm can solve complex non-differentiable optimization problems, which makes no assumptions regarding the targeted problem and can efficiently search a large feasible region [24, 38].

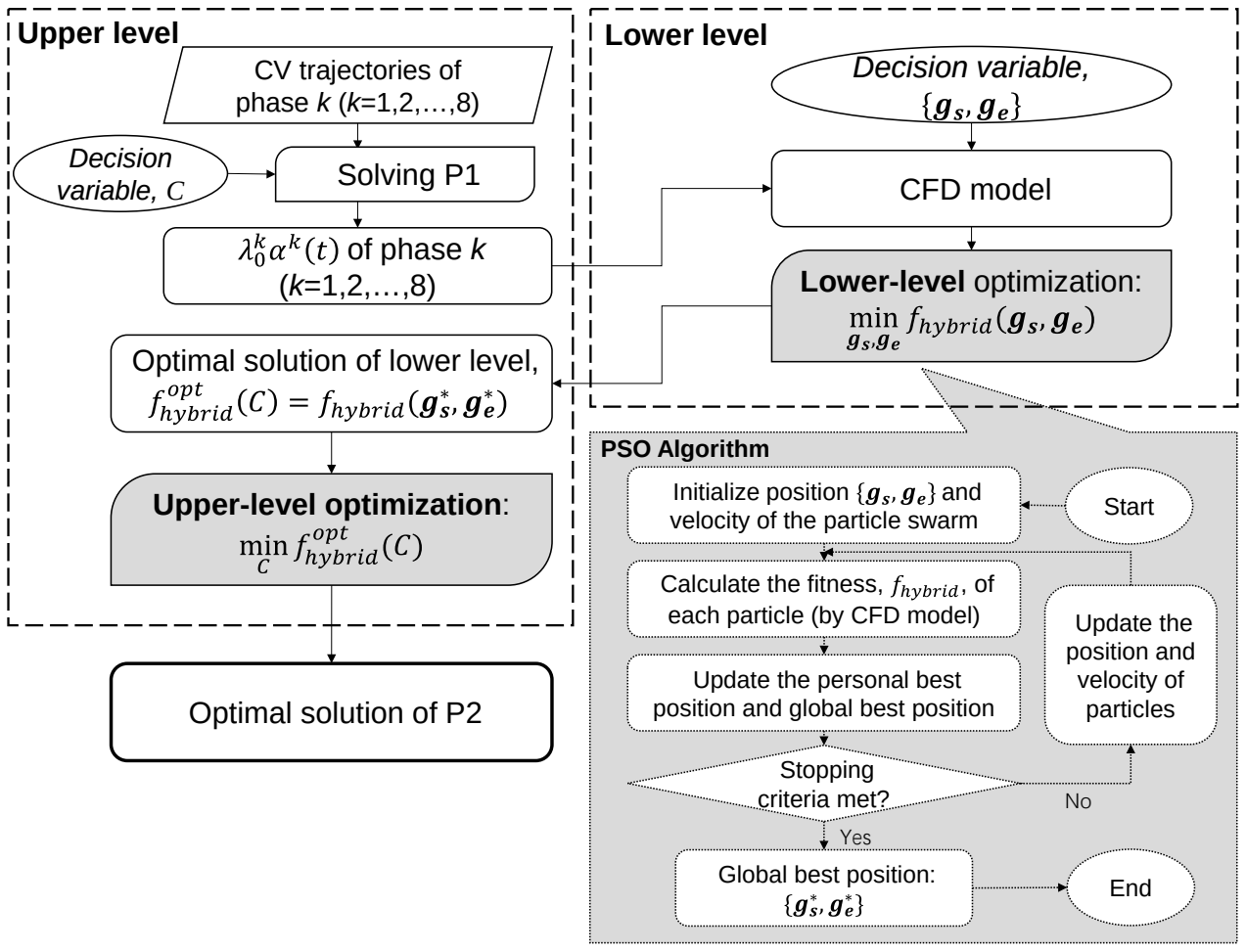

**Fig. 9.** Bi-level PSO-based algorithm to solve P2.

In the bi-level PSO-based algorithm, the cycle length and green splits are optimized separately, which allows us to less repeat the process of fitting $\alpha(t)$ by solving P1, thereby effectively improving the computational efficiency of the solution algorithm.

The objective functions for the two levels are as follows:

$$Upper\ level: \min_C f_{hybrid}^{opt}(C) \tag{34a}$$

$$Lower\ level: f_{hybrid}^{opt}(C) = \min_{\boldsymbol{g_s}, \boldsymbol{g_e}} f_{hybrid}(\boldsymbol{g_s}, \boldsymbol{g_e}) \tag{34b}$$

The process of the proposed solution algorithm is as follows: In the upper level, the CV trajectories are used as input. Given the decision variable *C*, the $\alpha(t)$ for all phases is fitted by solving P1 and then passed to the lower level. Then, at the lower level, the cycle length and $\alpha(t)$ for all phases are fixed. The PSO algorithm is applied to attain the optimal green splits and the corresponding minimum objective value, where the proposed CFD model is adopted to calculate the fitness of the particles. After obtaining the optimal solution of the lower-level optimization, that is, $f_{hybrid}^{opt}(C)$, for each feasible cycle length, the final optimal solution of P2 can be easily obtained in the upper-level optimization.

## IV. Evaluation

## *A. Test scenarios*

The performance of the proposed method was evaluated based on VISSIM simulation. To verify the effectiveness of the proposed method under both random arrival and platoon arrival patterns, a simulation model was built based on three consecutive intersections on Jinling Road, Changzhou City,

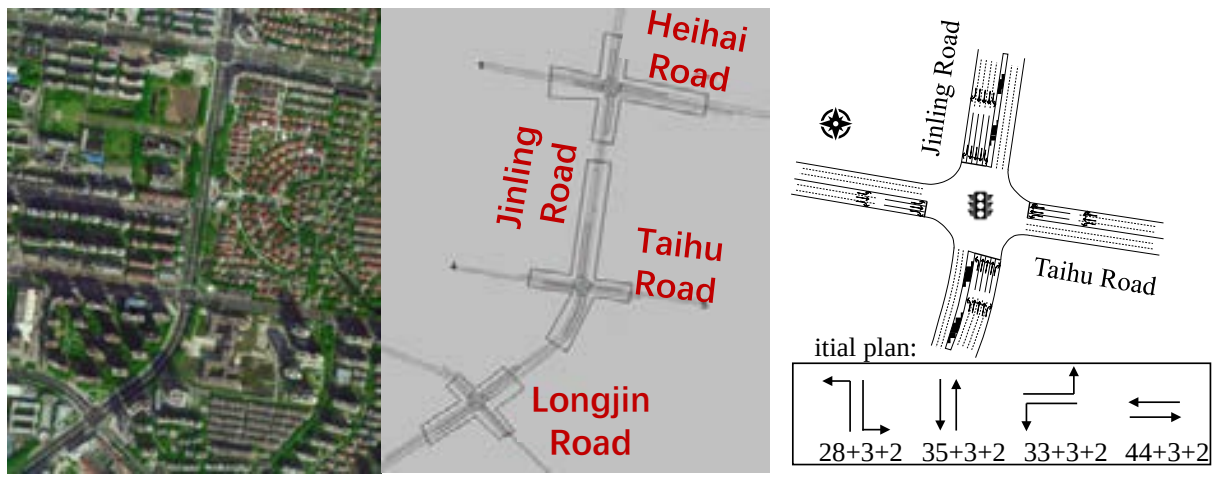

(a) Simulation model and intersection layout

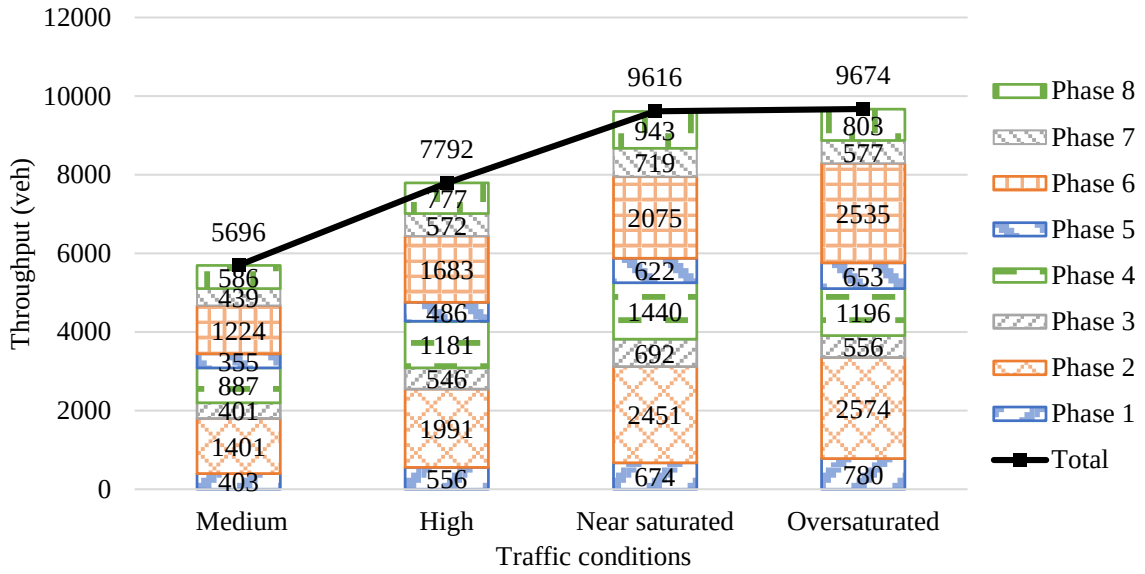

(b) Throughput of different traffic conditions

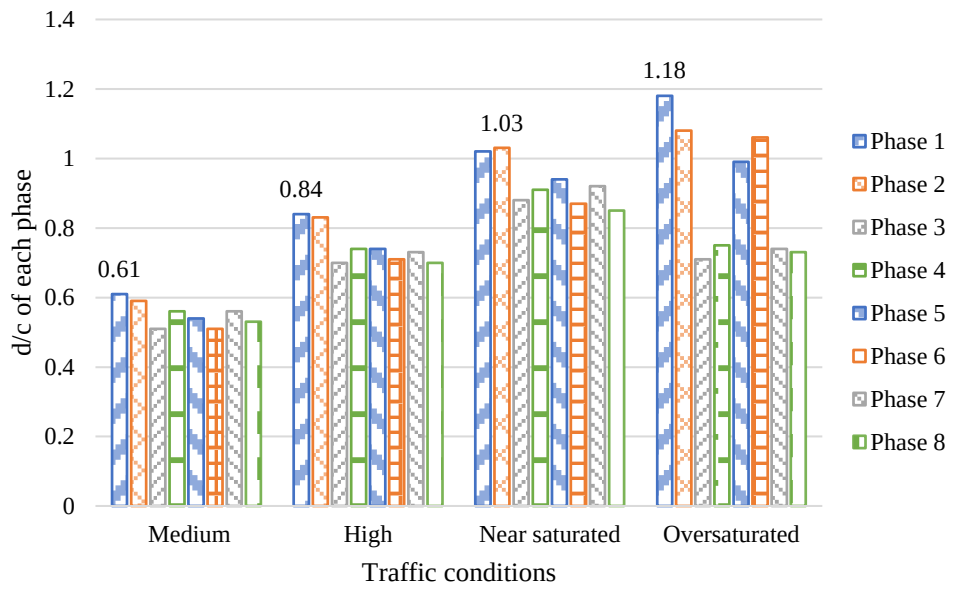

(c) d/c of different traffic conditions

**Fig. 10.** Basic information of the simulation model.

China, as shown in Fig. 10(a). The subject intersection is in the middle, that is, the Jinling-Taihu intersection. The lane allocation and the initial plan of the subject intersection are also shown in Fig. 10(a), which are provided by the local traffic management center. In this case, the right-turning flow was not controlled by signals. To pursue more realistic scenarios, the truck ratio and turning ratio were collected by the license plate recognition data deployed at the stop lines of three intersections and processed to calibrate the simulation model. The simulation model was run for 9000 s, where the first 900 s was to warm up and the last 900 s was to ensure the integrity of vehicle trajectories. Thus, two hours of data from 900 to 8100 s in the simulation were used for evaluation.

To better demonstrate the advantages of the proposed method, signal timing plans generated by Synchro were also applied in VISSIM and used for comparison. The lost time for each phase was 6 s, consisting of 3 s of start-up lost time, 1 s of an unused portion of the yellow time, and 2 s of red time. Four typical traffic conditions indicated by the demand to capacity (d/c) under the initial plan were tested in the simulation, as shown in Fig. 10(b) and (c), including medium (d/c = 0.6), high (d/c = 0.8), near-saturated (d/c = 1.0), and oversaturated (d/c = 1.2) conditions. The CVs were sampled from the population under different penetration rates, that is, 5–50%. Each penetration rate was randomly sampled 10 times with different sampling seeds to ensure result reliability. To better reflect the adaptability of the model to traffic flows, the time required for pedestrians to pass through the intersection was not considered in our cases; thus, the minimum green time was set as 5 s for all phases [7]. Notably, the input of Synchro is the ground truth volume output by VISSIM, whereas the input of the proposed method is the CV trajectories. All optimized plans are applied in VISSIM again, with all other settings maintained the same as those in the initial plan. The jam spacing $d_0$ is empirically set to 7 m in the evaluation.

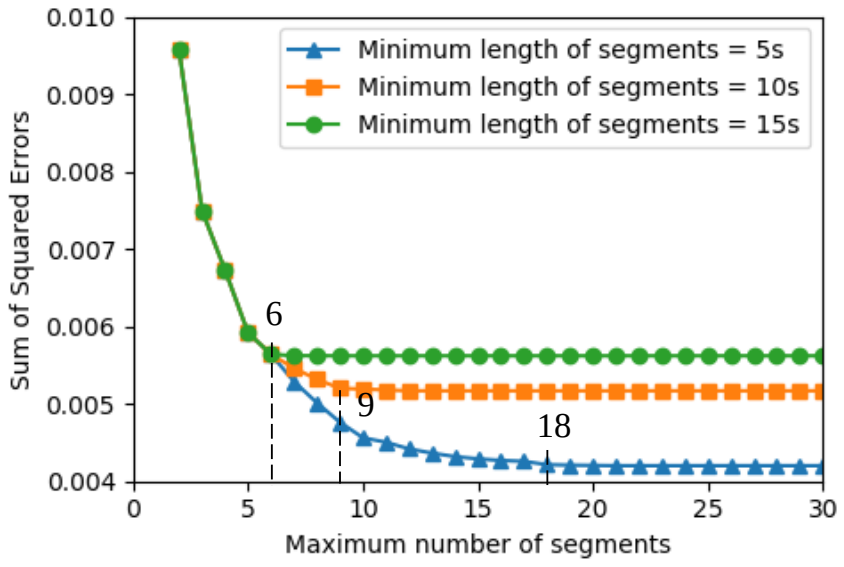

(a) SSE under different parameters

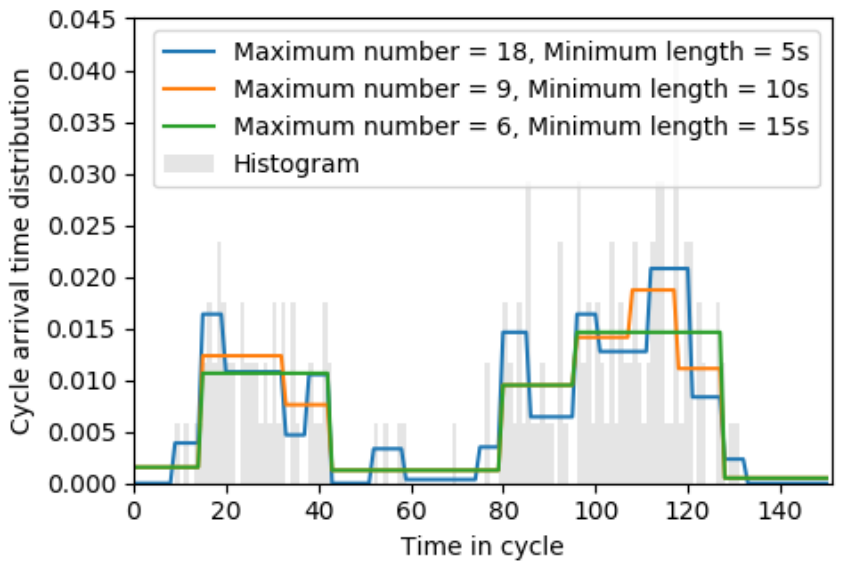

(b) Approximated $\alpha(t)$

**Fig. 11.** Test results of $\alpha(t)$ approximation.

### *B. Estimation of* $\lambda(t)$

Before we apply the optimization model, the estimation of $\lambda(t)$ must be evaluated, which includes the approximation of $\alpha(t)$ and estimation of the average arrival rate $\lambda_0$.

*1) Approximation of* $\alpha(t)$*:*

$\alpha(t)$ must be approximated based on CV trajectories by the regression tree method for each feasible cycle length, in which two critical parameters need to be predetermined, that is, the maximum number of segments and the minimum length of segments. This section aims to determine the most suitable combination of the two critical parameters.

We tested the approximation of $\alpha(t)$ based on a CV dataset collected from the real world at the Fuzhong-Huanggang intersection in Shenzhen, China. This dataset was used for the CFD estimation in Tan *et al.* [33]. The approximation results of $\alpha(t)$ for different combinations of parameters are shown in Fig. 11(a). By fixing the minimum lengths of segments, the SSE descends faster when the maximum number of segments is smaller. With an increase in the maximum number of segments, the SSE decreases while the descending rate becomes smaller until it remains unchanged. Thus, for different minimum lengths of segments, we can obtain the optimal maximum numbers of segments, which are 6, 9, and 18 when the minimum lengths of segments are 5 s, 10 s, and 15 s, respectively.

The approximation results of the three optimal combinations of the two parameters are shown in Fig. 11(b). As can be seen, with a smaller minimum length of segments, approximated $\alpha(t)$ can better accommodate the distribution of CV arrivals. However, notably, because the CV trajectories are randomly sampled from the population, the arrival distribution of CVs may not perfectly represent the arrival distribution of the population, especially when the penetration rates are low. In addition, considering the impact of the upstream intersection, vehicle arrivals at downstream intersections typically exhibit bimodal characteristics. Therefore, adopting a small minimum length of segments may introduce errors caused by the random sampling of CVs. In summary, the recommended and adopted combination of parameters for $\alpha(t)$ approximation in the proposed method is 9 for the maximum number of segments and 10 s for the minimum length of segments, which not only captures the characteristics of the CV arrival distribution well but is also less sensitive to fluctuations in the distribution.

*2) Estimation of $\lambda_0$:*

The average arrival rate $\lambda_0$ is a critical input of the CFD; Tan *et al.* [33] found that its estimation accuracy determines the accuracy of the CFD estimation and prediction. In other words, it can determine the reliability of the proposed CFD model. In this section, the estimation of $\lambda_0$ for all phases at the studied intersection is evaluated under different d/c and penetration rates.

The mean absolute percentage error (MAPE) was adopted to evaluate the accuracy of the estimates:

$$MAPE = \frac{1}{8}\sum_{k=1}^{8}\left|\frac{\hat{x}_k - x_k}{x_k}\right| \tag{35}$$

where $x_k$ is the ground truth value of phase *k* and $\hat{x}_k$ is the estimate of phase *k*.

The arrival rate estimation results are shown in Fig. 12. The line that indicates the performance after the distribution is obtained from Eq. (26), which uses the number of CVs to correct the initial estimates, that is, the line before distribution. Note that the redistribution process did not change the sum of the estimates, that is, the total volume of the intersection. As shown, for all values of d/c, the *MAPE* of the two estimates decrease as the penetration rates increase, particularly when the penetration rate is less than 15%. The accuracy was significantly improved after redistribution based on the number of CVs. In addition, the error bar, that is, the standard deviation of 10 estimation errors, shows that the estimates become more stable after redistribution. When the penetration rates are as low as 5%, the average *MAPE* for the 10 experiments is no more than 12% for 0.6 d/c and no more than 8% for a higher d/c. When the penetration rates are 10%, which is close to the real-world condition, the average *MAPE* is no more than 8% for 0.6 d/c and approximately 6% for a higher d/c.

In Fig. 13, an example showing the estimation error of each phase is provided when d/c is 1.2, and the penetration rate is 10%. The percentage error of $\lambda_0$ estimation of the eight phases shows that before redistribution, the maximum overestimation reaches 22.5% and the minimum underestimation reaches -12.5%. After redistribution, the range from underestimation to overestimation is reduced by 10.6%, from -7.9% to 16.5%. The results show that the estimation errors of most phases (6/8) are reduced after redistribution, leading to a smaller *MAPE* of the intersection on average. Thus, we can conclude that considering the number of CVs of phases for traffic volume estimation can significantly neutralize the estimation errors and improve overall accuracy.

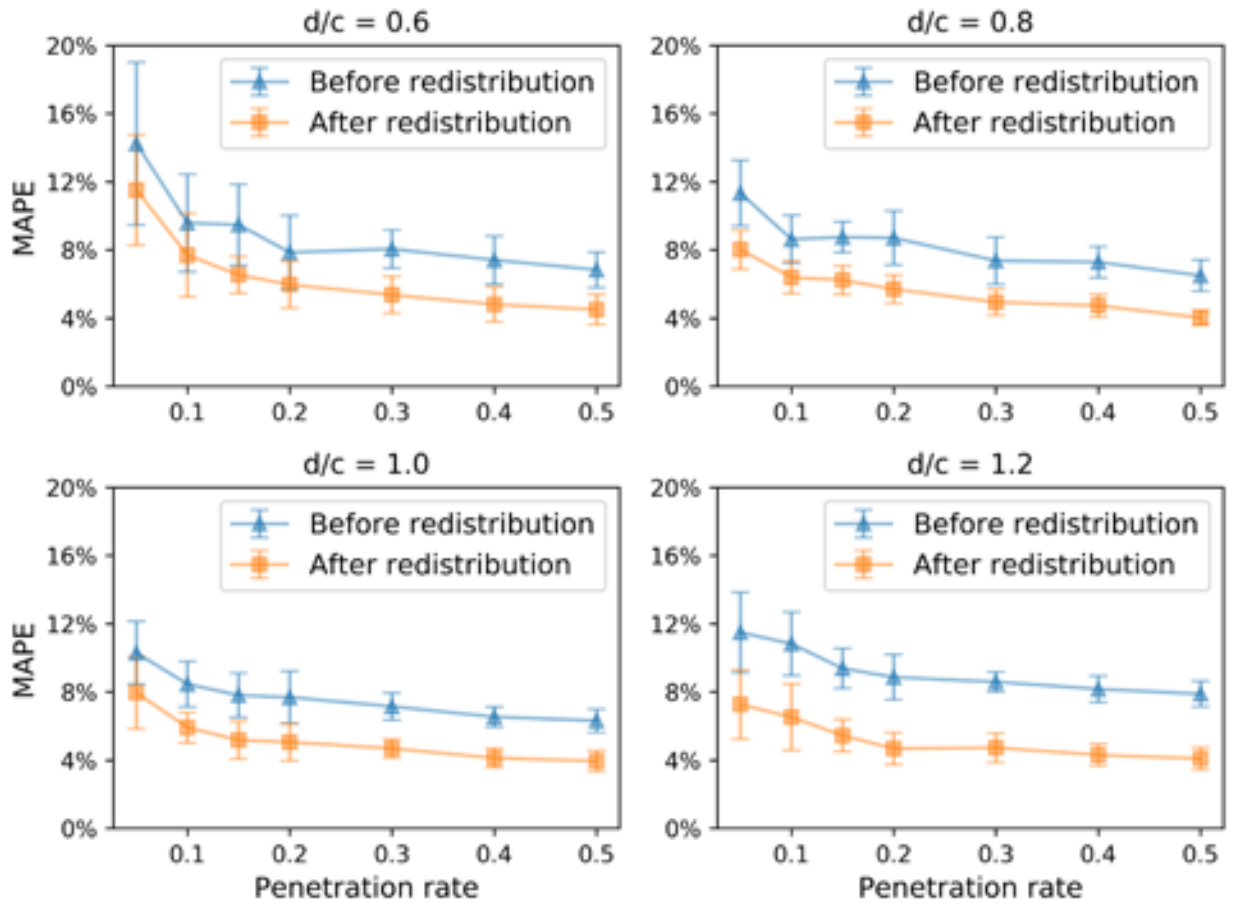

**Fig. 12.** Arrival rate estimation results.

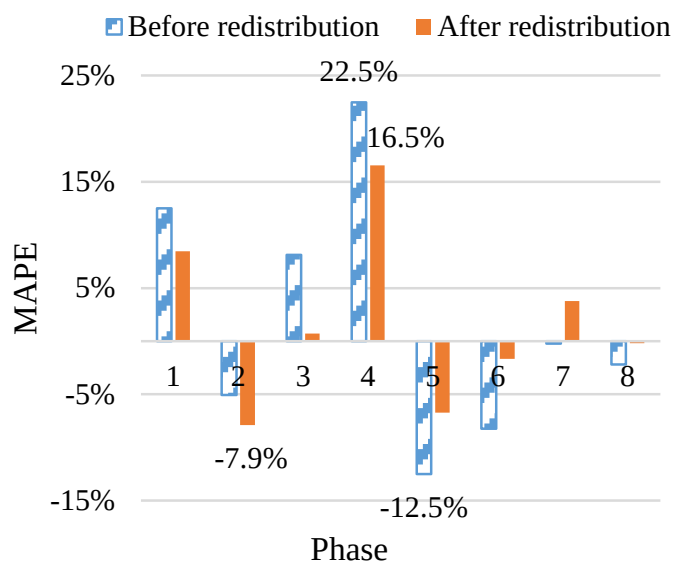

**Fig. 13.** Example of $\lambda_0$ estimation errors before and after redistribution (penetration rate = 10%, d/c = 1.2).

## C. Optimization results

After obtaining the average arrival rates for all phases, the CFD model can be applied to optimize the signal timing plan. As mentioned previously, the proposed optimization model was solved using a bi-level PSO-based algorithm. The upper level compared the optimal solution of the lower-level optimization for all feasible cycle lengths, while the lower level adopted the PSO method to attain optimal green splits. All experiments were conducted on a laptop with a 1.8 GHz 8-core i7-8550U CPU and 16 GB RAM. The model required approximately 30 s to attain the optimal green splits for the PSO method. The calculation speed of PSO was determined by the swarm size and stopping criterion. In our cases, the swarm size was set to 120, and the stopping criterion was that the optimal result remained unchanged for more than 50 iterations. Normally, the algorithm stops after approximately 70 iterations.

*1) Different d/c:*

First, the proposed optimization model was tested under different degrees of saturation. The penetration rate of the CVs, in this case, was 10%, which is close to real-world conditions. The improvement was used to evaluate the performance of the optimized plans compared with the initial plan.

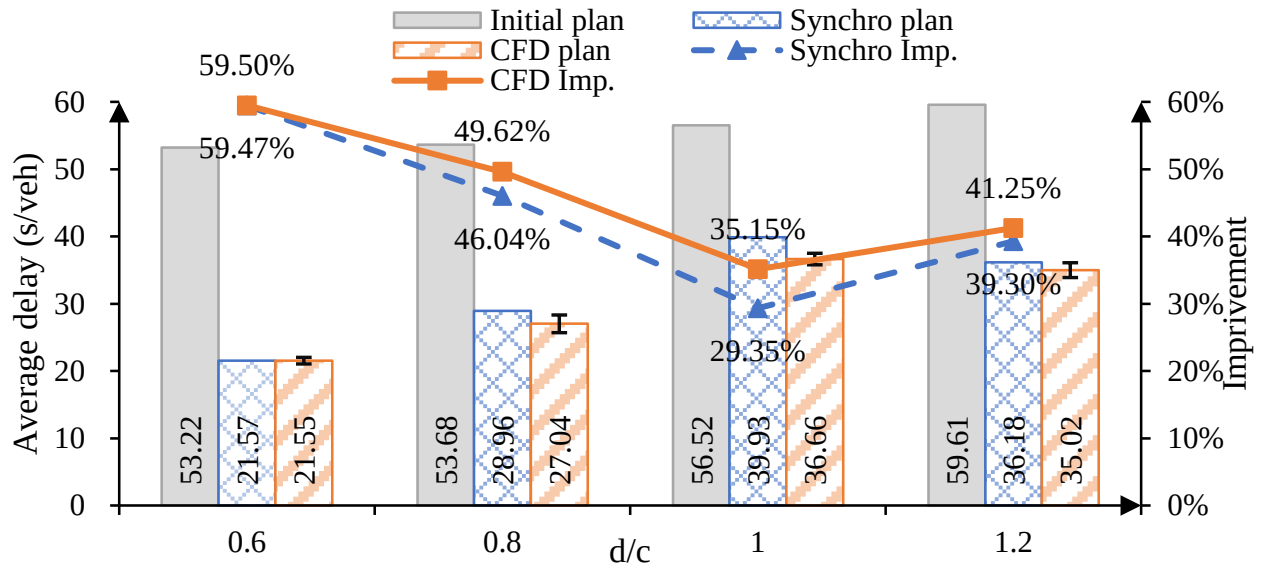

**Fig. 14.** Comparison of the average delay.

TABLE I
COMPARISON OF OTHER PERFORMANCE MEASURES

| Performance measures | d/c | Signal timing plan | | | *IMP* | |
|---|---|---|---|---|---|---|
| | | Initial | Synchro | CFD | Synchro | CFD |
| Average queue (m) | 0.6 | 18.98 | 7.81 | 7.63 | 58.85% | 59.79% |
| | 0.8 | 25.94 | 17.60 | 13.36 | 32.16% | 48.51% |
| | 1.0 | 35.45 | 26.60 | 23.46 | 24.96% | 33.84% |
| | 1.2 | 37.81 | 24.19 | 22.78 | 36.01% | 39.75% |
| Stops | 0.6 | 0.78 | 0.77 | 0.81 | 1.61% | -3.85% |
| | 0.8 | 0.80 | 0.76 | 0.79 | 5.51% | 1.61% |
| | 1.0 | 0.88 | 0.83 | 0.89 | 5.90% | -0.85% |
| | 1.2 | 0.95 | 0.84 | 0.96 | 11.78% | -0.48% |
| Throughput (veh) | 0.6 | 7617 | 7618 | 7613 | -0.01% | 0.05% |
| | 0.8 | 10467 | 10478 | 10474 | -0.11% | -0.07% |
| | 1.0 | 12866 | 12790 | 12818 | 0.59% | 0.37% |
| | 1.2 | 12820 | 12866 | 12864 | -0.36% | -0.34% |

$$IMP = \frac{x_{ini} - x_{opt}}{x_{ini}} \quad (36)$$

where $x_{ini}$ is the value of the performance measures under the initial plan, and $x_{opt}$ is the value under the optimized plans.

Fig. 14 shows the average delay of the initial plan, Synchro plan, and plan using the proposed model (denoted as the CFD plan). The results of the CFD plan were an average of 10 experiments. The error bars indicate that the proposed method can produce stable results for different CV samples. As shown, for all d/c, both the Synchro plan and the CFD plan can significantly improve the average delay compared with the initial plan, and the CFD plan performs better than the Synchro plan. When the d/c is 0.6, the average delays of the Synchro plan and that of the CFD plan are almost the same; the improvement of the CFD plan reaches 59.50%, which is only 0.03% higher than that of the Synchro plan. When d/c is 1.0, the delay improvement of the CFD plan is 35.15%, which is the smallest among all d/c. However, the delay improvement of the Synchro plan is even 5.8% less than that of the CFD plan, which is the largest difference among all d/c. When the d/c is 0.8, the delay improvements of the CFD plan and the Synchro plan are 49.62% and 46.04%, respectively. While for oversaturated conditions, i.e., when d/c is equal to 1.2, the delay improvements of the CFD plan and the Synchro plan are 41.25% and 39.30%, respectively.

Some other commonly used performance measures output by VISSIM, including the average queue, stops, and throughput of the CFD plan and the Synchro plan are shown in Table 1. Similar to the average delay, the improvement of the average queue of the CFD plan is also significantly greater than the Synchro plan for all d/c. As for stops, the Synchro plan performs better than the CFD plan for all d/c. Except for d/c = 0.8, the average number of stops of the CFD plan under other d/c increases compared with the initial plan. Note that the stops outputted by VISSIM also include stops caused by lane-changing behaviors, which are actually greater than stops caused by signal control. When the intersection is oversaturated with a d/c of 1.2, the throughputs of both the Synchro plan and CFD plan are increased.

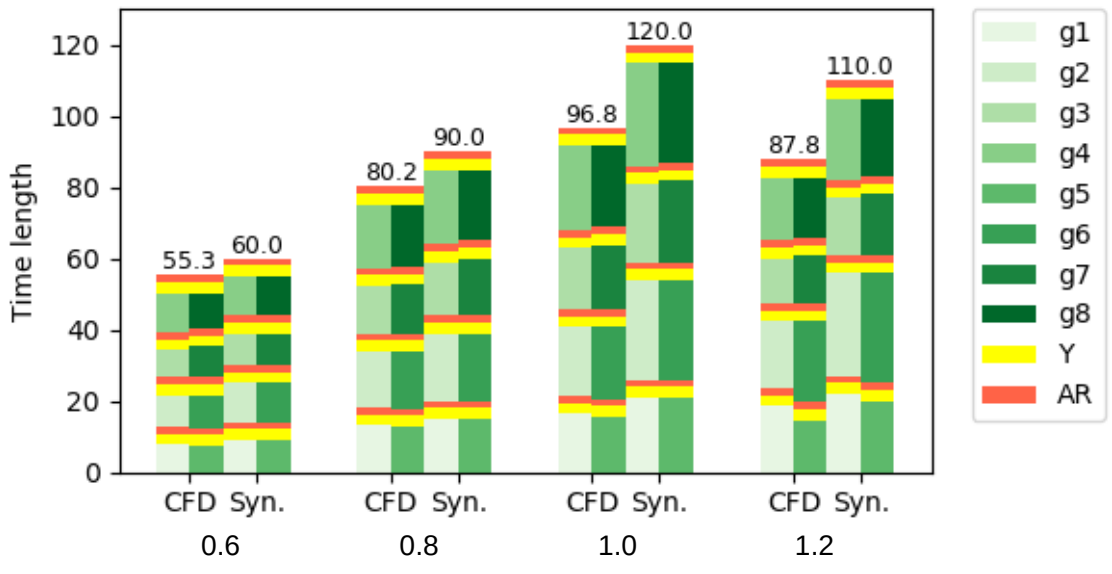

(a) Average plan

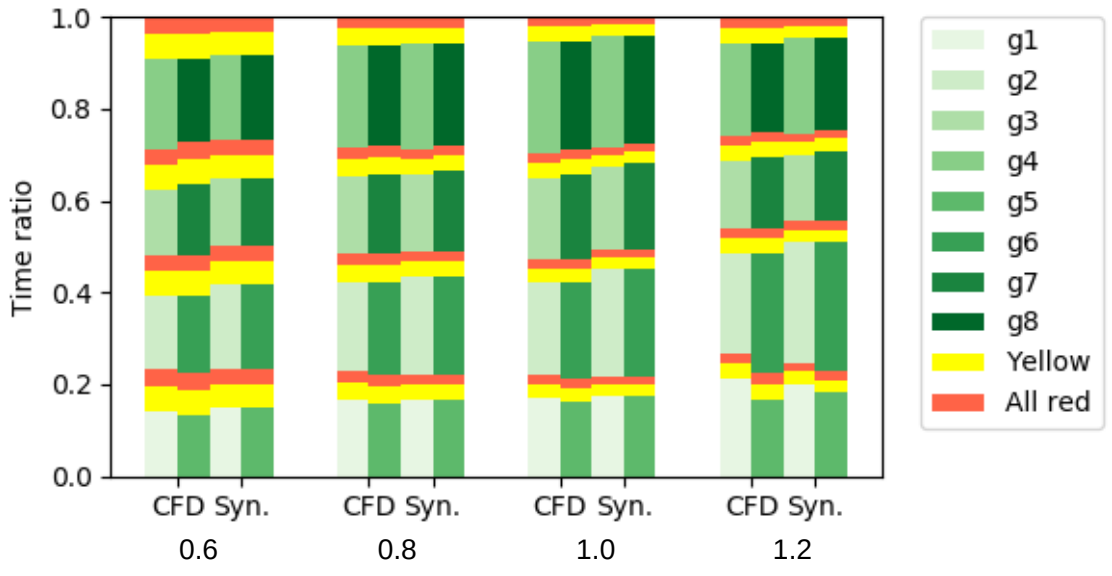

(b) Average ratio

**Fig. 15.** Comparison of the optimized plans.

The Synchro plan and the average CFD plan for the 10 experiments are shown in Fig. 15. The cycle length of the CFD plan is always smaller than that of the Synchro plan, and this difference is particularly significant at higher d/c. The green ratio shows that the CFD plan tends to allocate a greater green ratio to the second barrier in this case compared to the Synchro plan.

Fig. 16 shows a typical case of aggregated trajectories of a CFD plan and the Synchro plan when the d/c is 0.8 and the penetration rate is 10%. Trajectories in blue are vehicles that are not queued, in gray are vehicles queued once, and in red are vehicles queued more than once. The red and green points indicate the critical trajectory points when the vehicles joined the queue and left the queue, respectively. Vehicle arrivals of phases 1 are platoon arrivals affected by the upstream intersection, whereas vehicle arrivals in phase 8 are random arrivals. The cycle length of the CFD plan is 80 s, where that in the Synchro plan is 90 s.

As shown, owing to traffic demand fluctuations, the CFD plan and the Synchro plan both have a small proportion of vehicles queued twice. Besides, the cycle length of the CFD plan is smaller than that of the Synchro plan, and although the proportions of vehicles queued once and queued twice of phases 6 and 8 of the CFD plan are slightly greater than in the Synchro plan, the values of the average delay are still smaller. Notably, the phase 1 of the CFD plan show specific arrival patterns during the cycle because the 80 s of the cycle length is exactly half of the upstream intersection. The proportions of vehicles queued once and queued twice of the CFD plan are smaller than those of vehicles queued once and queued twice of the Synchro plan, and so is the average delay. This means that phase 1 received extra benefits from its fixed arrival patterns during the cycle. This also indicates that we can further improve the performance of the signal control at urban intersections by optimizing the reference point of the targeted intersection, which determines the arrival patterns during the cycle along with the cycle length.

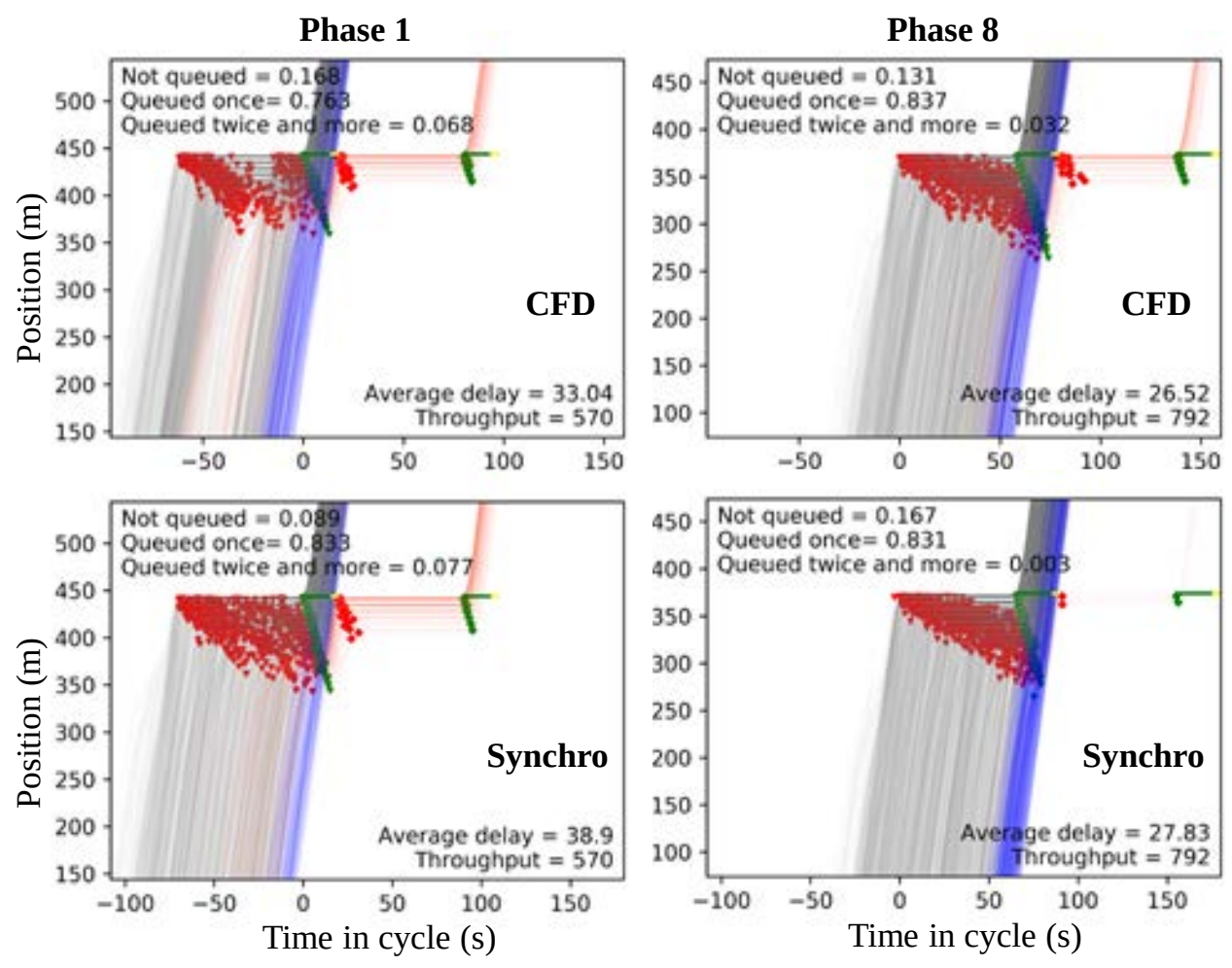

**Fig. 16.** An example of aggregated trajectories.

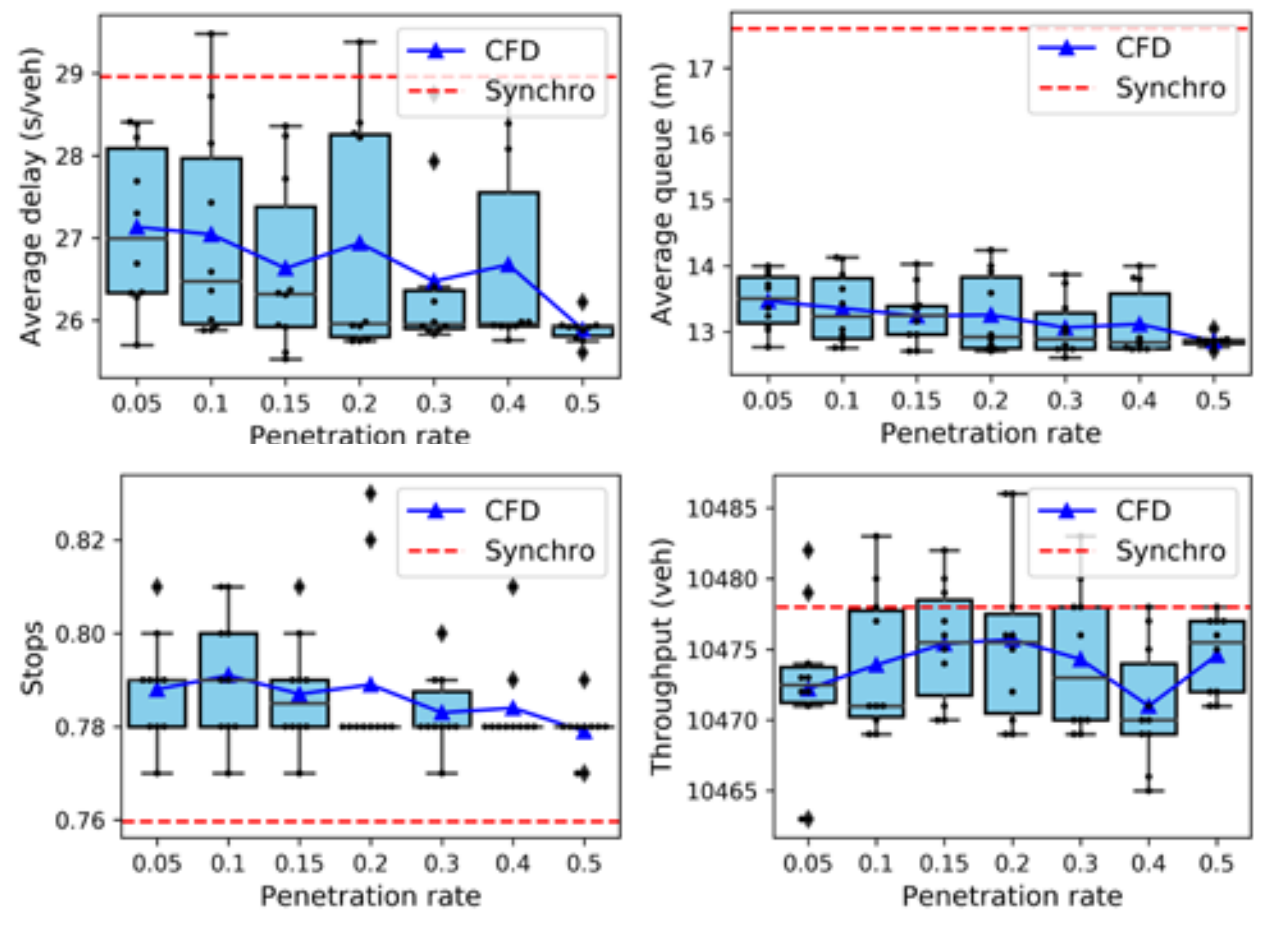

**Fig. 17.** Performance measures of the CFD plans under different penetration rates.

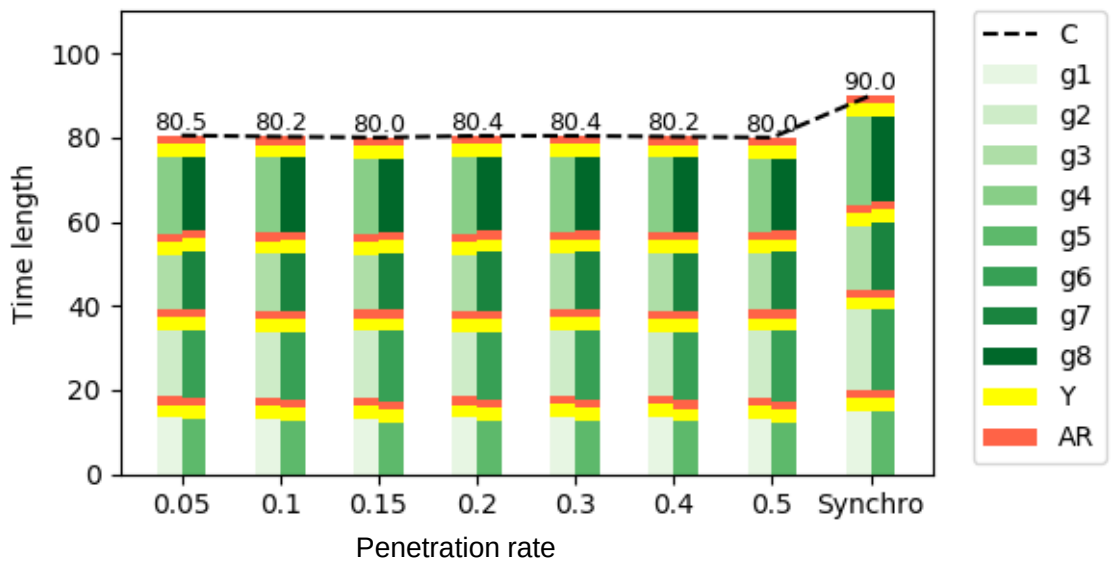

**Fig. 18.** Average CFD plans under different penetration rates.

*2) Different penetration rates:*

In this section, we evaluate the performance of the proposed optimization model under different penetration rates, and d/c was set to 0.8. As shown in Fig. 17, the overall trends of the average delay, average queue, and stops of the CFD plan decrease with increasing penetration rates, despite some fluctuations due to randomness of CVs among the 10 experiments. One exception is the throughput, which peaked at a 20% penetration rate and reached troughs at 5% and 40% penetration rates. Compared with the Synchro plan, the average delay and average queue of the CFD plan are significantly smaller. Nevertheless, the CFD plan performs slightly worse in terms of the stops and throughput. Notably, the variance of the average performance measures under different penetration rates was small. For instance, the average delay was 25.9 s/veh at the 50% penetration rate, while it only increased by 1.2 s/veh, i.e., 27.1 s/veh, at the 5% penetration rate.

The average CFD plans for the different penetration rates are shown in Fig. 18. The cycle lengths of the CFD plan under different penetration rates remained almost unchanged, at about 80 s, which is smaller than that of the Synchro, whose cycle length is 90 s. As is shown, the green splits of the CFD plan are also stable under different penetration rates.

TABLE II

COMPARISON OF PERFORMANCE MEASURES AFTER CONSIDERING THE REFERENCE POINT

| Signal timing plan | Average delay (s/veh) | Average queue (m) | Stops | Throughput (veh) |
|---|---|---|---|---|
| CFD0 ($C = 80, \Delta\varphi = 0$) | 25.88 | 12.88 | 0.78 | 10471 |
| CFD75 ($C = 80, \Delta\varphi = 75$) | 25.18 | 12.72 | 0.76 | 10472 |
| Syn0 ($C = 90, \Delta\varphi = 0$) | 28.96 | 17.60 | 0.76 | 10478 |
| Syn60 ($C = 80, \Delta\varphi = 60$) | 26.11 | 13.68 | 0.74 | 10475 |

*3) Considering the reference point:*

A noteworthy advantage of the proposed CFD-based optimization model is that the dynamic vehicle arrivals during the cycle are modeled by the CFD model based on CV trajectories. For urban intersections, vehicle arrivals at each entrance are, to some extent, determined by the signal of the upstream intersections. For coordinated intersections with a common cycle length, numerous studies have improved the performance of arterial coordination control by considering the offset during signal timing optimization. In this study, we believe that optimizing the relative relationship between the reference points of the targeted and upstream intersections would lead to extra benefits because the dynamic vehicle arrivals during the cycle are considered. This would still be the case even when the cycle length of the targeted intersection is different from that of the upstream intersections.

As indicated by Eq. (2), the evolution of vehicle arrivals during the cycle is associated with the cycle length and the

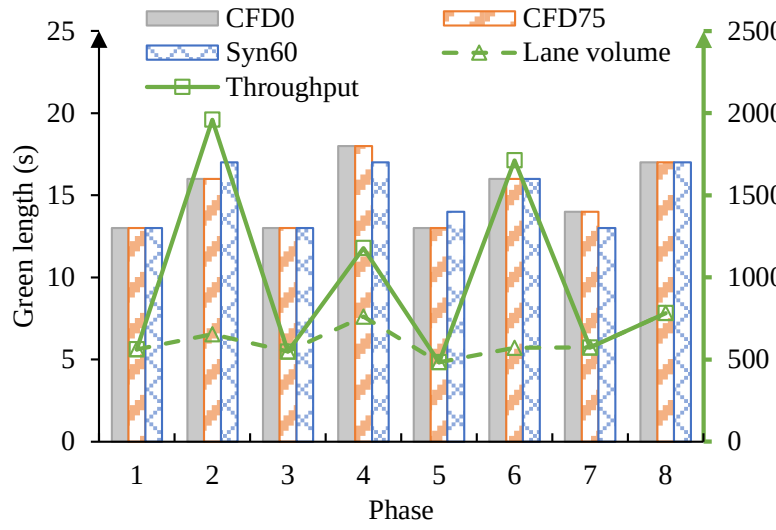

(a) Detailed plan

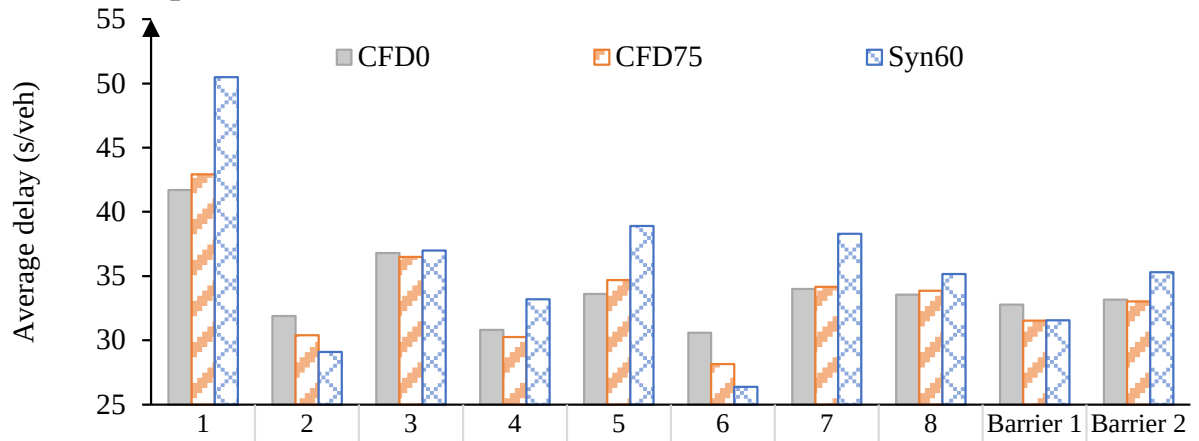

| | 1 | 2 | 3 | 4 | 5 | 6 | 7 | 8 | Barrier 1 | Barrier 2 |
|---|---|---|---|---|---|---|---|---|---|---|
| CFD0 | 41.71 | 31.89 | 36.81 | 30.8 | 33.62 | 30.59 | 34 | 33.55 | 32.78 | 33.16 |
| CFD75 | 42.92 | 30.39 | 36.5 | 30.27 | 34.7 | 28.15 | 34.16 | 33.86 | 31.54 | 33.02 |
| Syn60 | 50.49 | 29.09 | 36.99 | 33.19 | 38.9 | 26.37 | 38.29 | 35.15 | 31.57 | 35.31 |

(b) Average delay of each phase

**Fig. 19.** Comparison of average delay in terms of phases.

reference point of the signal timing plan of the analysis period. In the previous evaluation, the reference point was fixed at 0, that is, $\Delta\varphi = 0$, and only the cycle length $C$ was optimized in the upper level of the solution algorithm. In this section, both the cycle length $C$ and the reference point $\Delta\varphi$ were optimized in the upper level of the solution algorithm.

An example experiment with 0.8 d/c and 10% penetration rate was conducted. Note that when d/c was 0.8, the cycle length of the Synchro plan was 90 s. The Synchro software can only optimize the “offset” if the cycle length is half or the same as the cycle length of upstream and the downstream intersections (160 s). Thus, the Synchro plan, in this case, is generated with the cycle length fixed at 80 s, which was proved to be better than the case with a cycle length of 90 s by our previous results.

The performance measures of the optimal CFD plan without considering $\Delta\varphi$ (denoted as CFD0), CFD plan considering $\Delta\varphi$ (denoted as CFD75), Synchro plan without considering $\Delta\varphi$ (denoted as Syn0), and Synchro plan considering $\Delta\varphi$ (denoted as Syn60) are shown in Table 2. After we considered $\Delta\varphi$ in the optimization, the average delay, average queue, and stops of both CFD75 and Synchro60 decreased. However, the optimal $\Delta\varphi$ values of the CFD plan and the Synchro plan were different, which were 75 s and 60 s, respectively. The average delay and queue of CFD75 were significantly smaller than those of Synchro60, while the stops were slightly greater, and the throughput was slightly smaller.

To further evaluate the effect of considering $\Delta\varphi$ during the optimization, the detailed plans of CFD0, CFD75, and Synchro60 with the same cycle length and their average delay in terms of phases are shown in Fig. 19; thus, the effect of the cycle length on the results can be excluded. As shown, the green splits of CFD0 and CFD75 were the same, and Synchro60 was also similar. The green lengths of the four phases of Synchro60 are the same as those of CFD75, which are phases 1, 3, 6, and 8. The differences between the other phases were just 1 s.

In our NEMA phase structure, vehicle arrivals of phases in

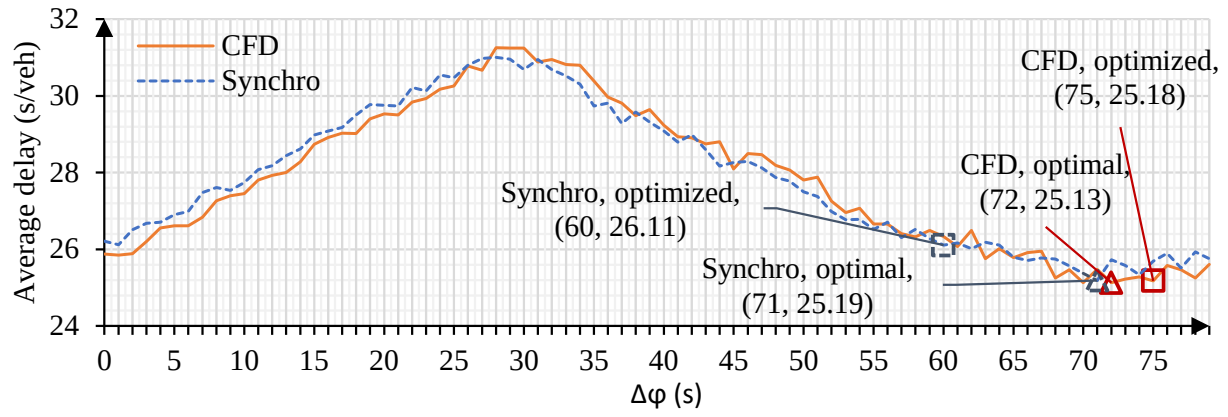

(a) Average delay under different $\Delta\varphi$ outputted by VISSIM

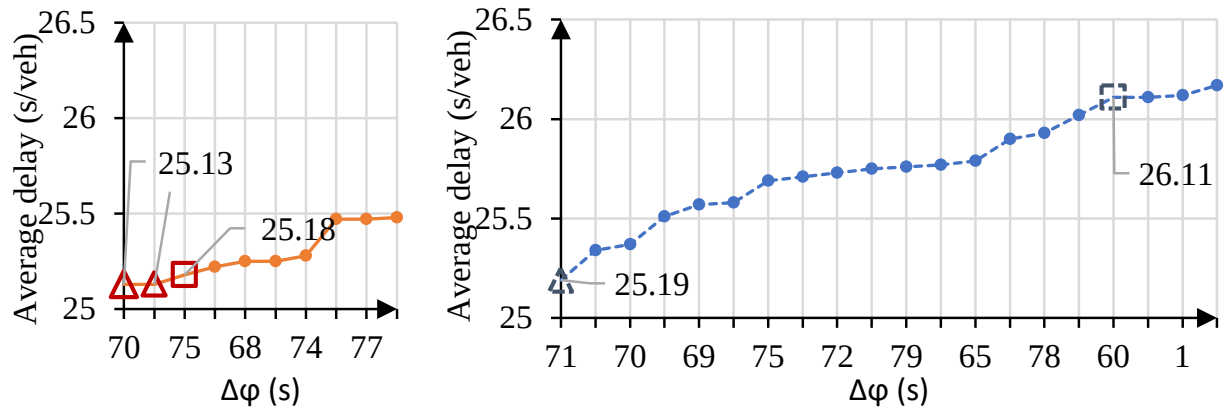

(b) Ordered $\Delta\varphi$ of CFD plans and Synchro plans

**Fig. 20.** Average delay of CFD plans and Synchro plans under different $\Delta\varphi$.

barrier 1 (i.e., phases 1, 2, 5, and 6) are affected by the signal control of upstream intersections and thereby are platoon arrivals, while the other phases (phases 3, 4, 7, and 8) in barrier 2 are random arrivals. Therefore, the effect of considering $\Delta\varphi$ must be discussed separately.

For phases with random arrivals in barrier 2, changing $\Delta\varphi$ has a limited effect on their performance, as indicated by the average delay of CFD0 and CFD75 in Fig. 19(b). The average delays of the second barrier of CFD0 and CFD75 are 33.16 s/veh and 33.02 s/veh, respectively. However, the average delay of the second barrier of Synchro60 is 35.31 s/veh, which is much greater. This is mainly caused by the smaller green lengths of phases 4 and 7 of Synchro60. For phases with platoon arrivals in barrier 1, we can see that the average delay of CFD75 is 31.54 s/veh, which is obviously smaller than that of CFD0, whose average delay of barrier 1 is 32.78 s/veh. Although Syn60 allocates a longer green length to barrier 1 than CFD75, its average delay of barrier 1 is still slightly greater than that of CFD75, which are 31.57 v/veh and 31.54 v/veh, respectively.

To validate whether the optimized $\Delta\varphi$ by the CFD-based model is optimal, the average delay of the CFD plan under different $\Delta\varphi$ is output by VISSIM, as shown in Fig. 20, where the cycle length and green splits are fixed. The triangle indicates the corresponding optimal $\Delta\varphi$ output by VISSIM, while the square indicates the $\Delta\varphi$ optimized by the proposed CFD-based model or Synchro. As shown, the optimal $\Delta\varphi$ for the CFD plan is 70 s or 72 s, whereas the optimized $\Delta\varphi$ is 75 s, which is the second-best value. However, $\Delta\varphi$ optimized by Synchro (60 s) is far from the corresponding optimal value (71 s). In summary, we can conclude that, for urban intersections where vehicle arrivals are affected by the signal control of upstream intersections, the proposed CFD-based optimization model can obtain an almost optimal $\Delta\varphi$, which is much better than that obtained by Synchro plans.

## V. Conclusion and Future Studies

In this study, we proposed a CFD-based signal timing optimization method for fixed-time signal control at isolated intersections, which is driven by the CV trajectory and does not require prior vehicle arrival assumptions. The CFD model profiles the vehicle arrival and departure processes under varying signal timing plans, in which the traffic demand is estimated based on the WMLE method. For signal timing optimization, a multi-objective optimization model was proposed for both undersaturated and oversaturated traffic conditions. The primary objective was to minimize the exceeded queue dissipation time of the intersection, aiming to dissipate the queue before the yellow start time as much as possible, while the secondary objective was to minimize the average delay of the intersection. Considering the data-driven property of the CFD model, a bi-level PSO-based algorithm was proposed for the solution. Notably, a significant advantage of the proposed method is that the dynamic vehicle arrivals during the cycle are considered and generated based on CV trajectories, which enables our method to not only attain better cycle length and green splits than those of Synchro but also further optimize the reference point of the signal timing plan, thereby achieving better performance.

The proposed method was evaluated based on simulation data and compared with Synchro, which uses the ground truth traffic volume as input. The results showed that the proposed method performed better than Synchro under varying traffic conditions in terms of the average delay and average queue length. As the penetration rates increased, the average delay, average queue length, and stops of the proposed method decreased. Even under a low penetration rate, i.e., 5%, the proposed method still outperformed Synchro in terms of the average delay and average queue length. When considering the reference point in optimization, both the average delay, average queue length, and stops of the CFD plan were further decreased, which also outperformed the Synchro plan.

In particular, we found that: 1) for arrival rate estimation, considering the number of CVs of phases can significantly neutralize the estimation errors and improve the overall accuracy; 2) the proposed CFD-based optimization model can obtain a smaller and better cycle length than those of Synchro under varying traffic conditions in terms of the average delay and average queue length; 3) under different penetration rates, the optimal plan by the proposed method is quite stable; 4) because the dynamic vehicle arrivals during the cycle are considered in the CFD model, some phases of the CFD plan can receive additional benefits from the fixed arrival patterns; 5) by optimizing the reference point, the performance of the proposed CFD-based optimization model can be further improved; and 6) the proposed CFD-based optimization model can obtain an almost optimal reference point, which is also better than that of the Synchro plan.

There are several possible directions for future research: 1) Currently, the proposed arrival rate estimation method simply used the number of CVs to redistribute the estimates and achieve a rough estimation of those phases without queued CVs, which is straightforward but rough. In the future, we will extend the current WMLE method to multiple phases and integrate the quantity information of CVs as prior information, and thus, develop a joint estimation method for multi-phase traffic demands based on Bayesian theory. 2) Although the CFD model can identify the oversaturated condition, it cannot accurately represent the periodic traffic flow profile of this traffic condition. Thus, future works can focus on the CFD model for oversaturated conditions. 3) Presently, the CFD model can only represent the average traffic flow profile while the traffic demand fluctuations are not considered. Therefore, future studies may improve the CFD model and consider cycle-to-cycle traffic demand fluctuations. 4) The current CFD model did not consider the influence of signal controls at neighboring intersections, but if we can further model the vehicle arrivals at downstream intersections based on the signal timing parameters of upstream intersections, the arterial and network signal control can be optimized. However, regarding the nonparametric nature of the CFD model, how to simplify the optimization model and solve it efficiently will be challenging.

## Credit Authorship Contribution Statement

Chaopeng Tan: Conceptualization, Methodology, Validation, Writing - original draft, Writing - review & editing. Yumin Cao: Validation, Writing - review & editing. Xuegang (Jeff) Ban: Conceptualization, Writing - review & editing, Supervision. Keshuang Tang: Conceptualization, Funding acquisition, Writing - review & editing, Supervision.

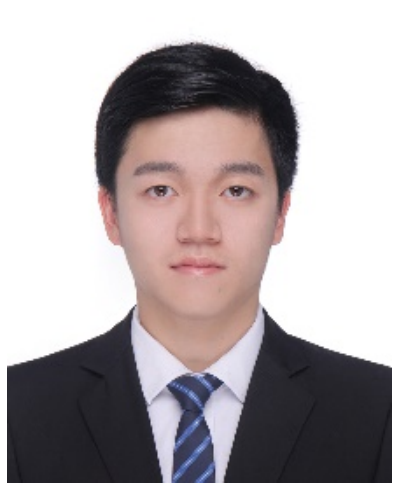

**Chaopeng Tan** received his Ph.D. degree in traffic engineering from Tongji University, Shanghai, China, in 2022. He was also a visiting Ph.D. student at the University of Washington (Seattle) for two years. His main research interests include urban traffic modeling and control with connected vehicles.

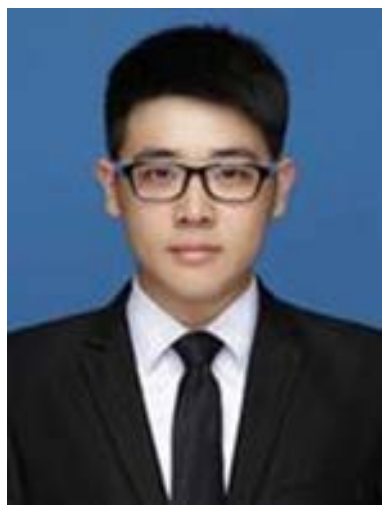

**Yumin Cao** received the B.S. degree in traffic equipment and control engineering from Central South University in 2018. He is currently pursuing the Ph.D. degree with the College of Transportation Engineering, Tongji University, China. His main research interests include path flow estimation and intelligent transportation systems.

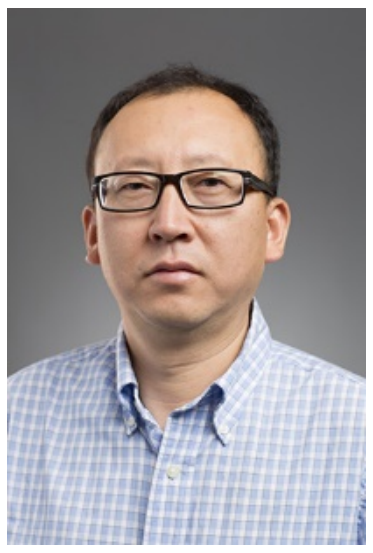

**Xuegang (Jeff) Ban** received the B.S. and M.S. degrees in automotive engineering from Tsinghua University and the M.S. degree in computer sciences and Ph.D. degree in transportation engineering from the University of Wisconsin–Madison. He is currently a Professor with the Department of Civil and Environmental Engineering, University of Washington. His research interests are in transportation network system modeling and simulation, urban traffic modeling and control, and intelligent transportation systems.

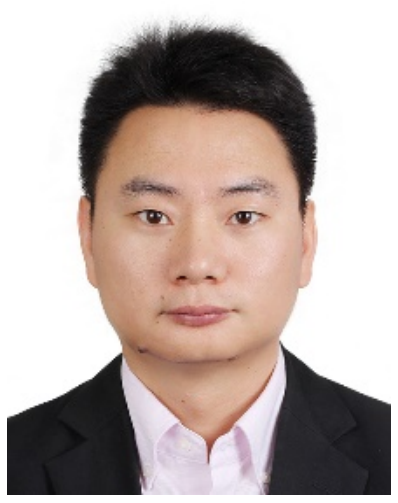

**Keshuang Tang** received his doctor’s degree in traffic engineering from Nagoya University in 2008. Afterward, he worked at The University of Tokyo as a postdoctoral research fellow, and then at Tohoku University as a Project Assistant Professor. He was also a visiting scholar of the University of California, Berkeley in 2010. He is currently a full professor at the College of Transportation Engineering, Tongji University, China. His main research interests include driver behavior, signal control, and Intelligent Transportation Systems.